\documentclass{amsart}

\usepackage{subcaption}
\renewcommand\thesubfigure{(\alph{subfigure})}
\usepackage{xcolor}
\usepackage[margin=1in]{geometry}
\usepackage{algorithm}
\usepackage{algorithmicx}
\usepackage{amssymb, mathrsfs,amsmath,graphicx, amsthm}
\graphicspath{{Figures/}}

\usepackage{booktabs}

\usepackage[noend]{algpseudocode}

\usepackage{tikz}
\usepackage{pgfplots}
\pgfplotsset{compat=1.18} 
\usepackage{pgfplotstable}
\usepackage{booktabs} 

\newtheorem{theorem}{Theorem}

\makeatletter
\def\algbackskip{\hskip-\ALG@thistlm}
\makeatother
\author{Christina Frederick}
\address{Department of Mathematical Sciences, New Jersey Institute of Technology, Newark, NJ 07102}
\email{christina.frederick@njit.edu}

\author{Haomin Zhou}
\address{School of Mathematics, Georgia Institute of Technology, Atlanta, GA 30332}
\email{hmzhou@math.gatech.edu}

\thanks{This work was supported by ONR N00014-24-1-2095 and NSF grants DMS-2307465 and DMS-2510829. CF also acknowledges the J. Tinsley Oden Faculty Fellowship Research Program.}

\title{ Multi-Agent Path-Planning in a Moving Medium via Wasserstein Hamiltonian Flow}

\begin{document}

\maketitle

\begin{abstract}
We present a finite dimensional variational model for multi-agent path-planning in which a group of agents traverses from initial positions to a target distribution in a moving medium. The model is derived using the agent-based formulation of the Wasserstein Hamiltonian flows that transport between probability distributions while optimizing a running cost. The objective is the mismatch between their final positions and the target distribution. The constraints are a system of Hamiltonian equations that provide the trajectories of the agents. The free variables on which the optimization is defined form a finite vector of the initial velocities for the agents. The model is solved numerically by the L-BFGS method in conjunction with a shooting strategy. Several simulation examples, including a time-dependent moving medium, are presented to illustrate the performance of the model.
\end{abstract}

\section{Introduction}

Multi-agent path-planning has been an extensively studied research topic in the past few decades. The main task is to find optimal trajectories for a multi-agent fleet from their initial locations to desirable target positions. Unlike common set-ups in the existing multi-agent path-planning, in this paper we consider a modified task, in which the initial and/or target positions of agents are given by distributions, and the agents traverse in a moving medium, such as an ocean current. 

\subsection{Problem Formulation}  \label{problem.formulation}

We consider a multi-agent system of $N$ agents, whose trajectories are denoted by $X_i(t) \in \mathbb{R}^d$, $1\leq i \leq N$, navigating in an unbounded planar environment characterized by a given background flow $w(t, x) \in C^1([0, T]\times \mathbb{R}^d; \mathbb{R}^d)$. We assume that the agents are holonomic without inter-agent communication used by the controller. The initial positions of the agents $X(0)$ are given and they are samples of a probability distribution characterized by a density function $\mu(x)$. The goal is to find trajectories of $X(t)$, $t\in [0,T]$, in the flow $w$ such that their final positions $X(T)$ follow a desirable distribution given by the density $\nu(x)$, while the paths optimize the running cost in expectation, such as the kinetic energy contributed by the agents. 

This modified task is different from commonly studied path-planning problems in several aspects. The most notable one is that the target locations are not given {\it a priori}. Instead, they are samples of the distribution $\nu$. Since both the initial and target locations are random samples, the trajectories $X_i(t)$ are consequentially random variables with density denoted by $\rho(t,x)$. Hence, the running cost is defined by the expectation with respect to $\rho$. Furthermore, background flow $w$ adds complications in determining the running cost. For example, the kinetic energy to be optimized should exclude the influence of $w$. 

It is not easy to adapt the methods developed for the commonly studied multi-agent systems to this new task. For example, a straightforward approach is to sample the final locations $X_i(T)$ according to $\nu(x)$ first, then applying the existing methods to find the paths for the sampled positions. However, many samples are needed to calculate the running cost, leading to expensive computations. Hence, those modifications in the task demand new models and methods to address the challenges faced by the existing strategies.









\subsection{Relevant literature}

The path-planning for multi-agent systems has been extensively studied in the past few decades. There are many classical methods such as A* \cite{LaValle2006PlanningAlgorithms,Hart1968APaths}, rapidly-exploring random tree (RRT) \cite{LaValle1998Rapidly-exploringPlanning}, and model predictive control (MPC) \cite{Kouvaritakis2016ModelControl}. They were originally developed for the single agent path-planning and have been extended to the multi-agent systems \cite{Bono2022AEnvironments}. Our study is closely related to some topics, such as density control, whose aim is to design control strategies that transport one density to another while minimizing desirable running costs \cite{Chen2016OptimalI,Chen2016OptimalII,Chen2018OptimalIII}, and mean field games \cite{Ruthotto2020AProblems}, in which the density function plays a vital role in the dynamics of a game with an infinite number of players. In recent years, the mean field game formulation has been used for path planning for a swarm of agents \cite{Kang2020JointGame, Onken2023AFinding}.   

Path-planning in a moving medium has a long history that can be traced to Zermelo's navigation problem \cite{Zermelo1931UberWindverteilung}. There are numerous underwater applications, e.g., control of marine robots \cite{Hong2025ControlIntelligence, Ouerghi2022ImprovedMapping}.
Reinforcement learning is another popular method advanced rapidly in past few years due to the development of machine learning \cite{Sartoretti2019PRIMAL:Learning, Tang2025DeepSuccesses}. They motivate us to consider this modified task, although we take a different approach in this study.

Our formulation is inspired by the recent mathematical developments in optimal transport and Wasserstein Hamiltonian flow (WHF)\cite{Chow2020WassersteinFlows}, in which a Hamiltonian system on the Wasserstein manifold, the probability density space equipped with the optimal transport distance, is derived as the critical point for the optimization problem with constraints. An agent-based interpretation of WHF has been introduced in \cite{Cui2022AEquation, Cui2022TimeFlows}. They provide theoretical support for our model, the derivation, and the algorithm. 
It is also worth mentioning that there have been several studies aiming to apply optimal transport formulation to robotics and multi-agent path planning \cite{Chen2024DensitySystems,  Chen2015OptimalSystem, Chen2021OptimalControl, Krishnan2018DistributedSwarms, LeAcceleratingTransport}.  








\subsection{Contributions and organization}

In this paper, we propose an agent-level optimization formulation to solve the path-planning problem described in Section \ref{problem.formulation}. We accomplish this goal by the following three steps. First, we modify the agent-based optimal transport formulation to account for the background flow $w$ when calculating the kinetic energy contributed by the agents. Second, we derive the corresponding Wasserstein Hamiltonian flow (WHF), which is a dynamic system to describe the trajectories of the agents. In the end, we establish a finite dimensional optimization with constraints, in which the objective is the mismatch of agent locations with the target distribution while the constraint is the derived WHF. The independent variables on which the optimization is defined are a collection of vectors for the initial velocities of agents, which are of finite dimension. The optimal trajectories are computed by the standard optimization methods such as L-BFGS, in conjunction with classical numerical schemes for ordinary differential equations (ODEs) like the 4th order Runge-Kutta (RK45) method. The main contributions include

\begin{enumerate}
    \item establishing the WHF in a moving medium; and 
    \item proposing a finite dimensional optimization with the WHF as constraints for the path-planning of the multi-agent system. 
\end{enumerate}

Compared with other path-planning methods, the proposed framework has several distinctions. The search space is the set of initial velocity vectors for the agents, not the entire trajectories, which are of infinite dimension, as commonly considered in the existing literature. This is a significant dimension reduction. The model is derived from the optimal transport and WHF. They provide the theoretical underpinning for the proposed model. The formulation is conceptualized at the density level, which is more suitable to account for the changes or uncertainties in the agent positions, while the computation is carried out at the agent level, which is more resource efficient.   

The paper is organized as follows. Section \ref{sec:whflows} provides the derivation of WHF in a moving medium. The agent based formulation is given in Section \ref{sec:agent}. The density and agent-based optimization formulations and the associated algorithm are presented in Section \ref{sec:optimization} and Section \ref{sec:numericalimplementation} respectively, and numerical experiments are given in Section \ref{sec:numerics}. We conclude with a number of potential extensions and applications of our model in Section \ref{sec:future}.






\section{Wasserstein-Hamiltonian flows in a moving medium}\label{sec:whflows}

The theoretical motivation is based on the dynamic formulation of optimal transport. The Wasserstein-2 distance between two probability measures $\mu$ and $\nu$ is equivalent to solving the variational problem known as the Benamou-Brenier formula \cite{Benamou2000AProblem}
\begin{align} \min_{v,\rho} \left\{\int |v|^2\rho(t,x)\,dx\, dt, \quad \frac{\partial \rho}{\partial t}+\text{div}(v\rho) =0,  
   \rho(0,x)=\mu(x), \rho(1,x)=\nu(x)\right\}. \label{eq:dynamicwass}
\end{align}
This setup models transport in a simple, homogeneous environment with no background current.

We first modify the continuity equation in \eqref{eq:dynamicwass} to account for a background current $w(t,x)$ in addition to the velocity field $v$, so the total flow that is moving the group of agents is $v+w$.

Suppose that $(v,\rho)$ is a pair satisfying 
\begin{align}
 &\quad \frac{\partial \rho}{\partial t}+\text{div}((v+w)\rho) =0,  \qquad    \rho(0,x)=\mu(x),\quad \rho(1,x)=\nu(x).\label{eq:contH}
\end{align}
The total velocity of the group of agents is given by $v_{\mathrm{total}}=v+w$ and the agent-imposed (control) velocity field is given by $v$. The control cost is the kinetic energy contributed by the agent-imposed velocity field $v$ 
\begin{align}
\mathcal{K}(v,\rho)=\frac{1}{2}\int |v(t,x)|^2\rho(t,x)\,dx\, dt. \label{eq:K}
\end{align}
We propose to solve the following minimization problem:
\begin{align}
    \mathcal{J}_F[\mu, \nu]=\min_{\rho, v} \mathcal{K}(v,\rho)  \quad \text{such that }\rho, v \text{ satisfy }\eqref{eq:contH}. \label{eq:JHdist}
    \end{align}

 Finding the optimizer of \eqref{eq:JHdist} is not a trivial task. It involves finding Wasserstein-type geodesic curves. In high dimensions, e.g., machine learning principles are used to accomplish this \cite{Gai2021ASpace,LiuLearningGeodesics}. Here, to solve \eqref{eq:JHdist}, we consider the Lagrange multiplier technique as in \cite{Cui2022AEquation}:
\begin{align*}
\begin{split}
 \mathcal{L}(\rho, v, \lambda) &=     \frac{1}{2} \int_{0}^1\int_{\mathbb{R}^d} |v(t,x)|^2\rho(t,x),\,dt\,dx \\
&+\int_{0}^1\int_{\mathbb{R}^d}\varphi(t,x)\left(\rho_t(t,x)+\nabla\cdot(\rho(t,x)(v(t,x)+w(t,x)))\right)\,dt\,dx.
    \end{split}
    \end{align*}
Here $\varphi$ is the Lagrange multiplier.

The critical point of $\mathcal{L}$ satisfies the Karush–Kuhn–Tucker (KKT) conditions:
\begin{align*}
\frac{\partial \varphi}{\partial t}(t,x)&=  -v(t,x) \cdot \left(\frac{1}{2} v(t,x)+w(t,x)\right) \\
\frac{\partial \rho}{\partial t}(t,x) &=-\nabla \cdot ((v(t,x)+w(t,x))\rho(t,x)) \\
v(t,x) &= \nabla \varphi(t,x). 
\end{align*}

 When $\rho(t,x)$ is positive for all $t$, the system can be viewed as a Wasserstein-Hamiltonian flow \cite{Chow2020WassersteinFlows} with the Hamiltonian map $\mathcal{H}(\rho,\varphi, t)$ defined by    \begin{align}\mathcal{H}(\rho,\varphi, t) =  \int_{\mathbb{R}^d}  \left(\frac{1}{2}|\nabla \varphi(t,x)|^2+ \nabla \varphi(t,x) \cdot w(t,x)\right) \rho(t,x) \, dx
\label{eq:EulerianHamiltonian}
   \end{align}
satisfying the Hamilton-Jacobi and continuity equations
   \begin{align}
\begin{cases}\dfrac{\partial \varphi}{\partial t}&=\dfrac{\partial\mathcal{H}}{\partial \rho}=-\nabla \varphi \cdot \left(\frac{1}{2}\nabla \varphi +w\right)\\
\dfrac{\partial \rho}{\partial t}&=-\dfrac{\partial\mathcal{H}}{\partial \varphi}=-\nabla \cdot ((\nabla \varphi +w)\rho).
\end{cases} \label{eq:optpdes}
\end{align}


 As a remark, if the background flow $w$ is time-independent, the Hamiltonian $\mathcal{H}$ defined in \eqref{eq:EulerianHamiltonian} is a constant. Indeed, along the optimal trajectory $(\rho_t, \varphi_t)$, we have:
\begin{equation*}
\frac{d}{dt} \mathcal{H}(\rho, \varphi) = \int_{\mathbb{R}^d} \frac{\delta \mathcal{H}}{\delta \rho} \frac{\partial \rho}{\partial t} + \frac{\delta \mathcal{H}}{\delta \varphi} \frac{\partial \varphi}{\partial t} \, dx = \int_{\mathbb{R}^d} \dfrac{\partial \varphi}{\partial t} (-\dfrac{\partial \rho}{\partial t}) + \dfrac{\partial \rho}{\partial t} (\dfrac{\partial \varphi}{\partial t}) \, dx = 0.
\end{equation*}
In the presence of a time-dependent flow $w(t, x)$, the Hamiltonian evolves according to the partial derivative of the background field:
\begin{equation*}
\frac{d}{dt} \mathcal{H} = \int_{\mathbb{R}^d} \rho(t,x) \left( \nabla \varphi(t,x) \cdot \frac{\partial w}{\partial t}(t,x) \right) dx.
\end{equation*}

Our approach in solving this modified optimal transport problem, described in the next section, is to find the initial $\varphi_0$ such that the trajectory of \eqref{eq:optpdes} starting at $(\mu, \varphi_0)$ passes through $\nu$ at $t=1$ where the terminal potential $\varphi(1,x)$ is implicitly determined by the target density. This is an analog of the well-known geodesic equation between the densities $\mu$ and $\nu$ on the Wasserstein manifold (see, e.g., \cite{Villani2003TopicsTransportation, Ambrosio2008GradientMeasures}).

As a remark, since $\varphi$ is defined up to an arbitrary constant, the system \eqref{eq:optpdes} does not have a unique solution, so we will formulate the problem in the $(\rho, v)$ variables as in \cite{Cui2022AEquation}.

\section{Agent-based formulation}\label{sec:agent}

This section reformulates the energy minimization problem for agent trajectories \eqref{eq:JHdist} into a finite-dimensional optimization problem. We convert to Lagrangian coordinates following the fluid dynamics formulation of Monge–Kantorovich optimal transport \cite{Benamou2000AProblem}. Assume $\mu,\nu$ are probability measures on $\mathbb{R}^d$ that are absolutely continuous with respect to Lebesgue measure, and that $\rho,v$ are solutions of the continuity equation with background flow $w$, given by \eqref{eq:contH}, that are sufficiently smooth (e.g. $v+w$ is Lipschitz in space and continuous in time).

To solve the minimization problem \eqref{eq:JHdist} from a particle-tracking perspective, we introduce the flow map $X: [0,1] \times \mathbb{R}^d \to \mathbb{R}^d$. Let $\xi\in \mathbb{R}^d$ be the initial position of an agent at $t=0$; its position at time $t$ is denoted by $X(t)$. In this framework, the total velocity of an agent is the sum of the agent-imposed control velocity $v$ and the background flow $w$:
\begin{equation*}
    \dot{X}(t) = v(t, X(t)) + w(t, X(t)), \quad X(0) = \xi.
\end{equation*}

If $\xi$ is obtained by sampling according to $\mu$, $X(t)$ is a random variable whose density is denoted by $\rho(t,x)$. By the conservation of mass, the density $\rho(t, \cdot)$ is the push-forward of the initial measure $\mu$ by map $X(t, \cdot)$. Substituting the control variable $v = \dot{X} - w$ into the cost functional \eqref{eq:K}, the problem becomes one of minimizing the Lagrangian action for each agent:
\begin{equation*}
    \mathcal{J}_w[\mu, \nu] = \min_{X} \int_{\mathbb{R}^d} \int_{0}^1 \frac{1}{2} \left| \dot{X}(t) - w(t, X(t)) \right|^2 dt \, d\mu(\xi),
\end{equation*}
subject to the endpoint constraints $X(0) = \xi$ and $X(1)_{\sharp} \mu = \nu$.

To derive the dynamics, we define the Hamiltonian $\boldsymbol{H}(x, v, t)$ via the Legendre transform of the Lagrangian $\boldsymbol{L}(x, \dot{x}, t) = \frac{1}{2}|\dot{x} - w(t, x)|^2$. Taking $v$ as the canonical momentum conjugate to $x$:
\begin{align*}
    \boldsymbol{H}(x, v, t) &= \sup_{\hat{v}} \left\{  v \cdot \hat{v} - \boldsymbol{L}(x, \hat{v}, t) \right\}\nonumber \\
    &= \frac{1}{2}|v|^2 + v\cdot w(t, x). \label{eq:LagrangianHamiltonian}
\end{align*}


The optimal trajectories $(x(t), v(t, x(t)))$ are governed by the Hamiltonian system. The following result can be derived from the standard theory of Hamiltonian mechanics, and a similar procedure is given in \cite{Cui2022AEquation}.

\begin{theorem}[Wasserstein-Hamiltonian Flow in a Moving Medium]
Let $\mu, \nu$ be probability measures on $\mathbb{R}^d$ with bounded second-order moment and let $w(t, x)$ be a $C^1$ vector field. The critical points $X(t)$ and $q(t)$ minimizing the control energy functional
\begin{eqnarray}
    \mathcal{A}(X(t), q(t)) = \mathbb{E}_{\mu}[\int_0^1  \frac{1}{2} |q(t)|^2 dt], \\
    \text{subject to }  \dot{X}(t) = q(t) + w(t, X(t)),\qquad  X(0)\sim \mu, \qquad X(1)_{\sharp} \mu = \nu
    \label{eq:action}
\end{eqnarray}
 are characterized by the Hamiltonian system:
\begin{equation}
    \begin{cases}
        \dot{X} =\frac{\partial \boldsymbol{H}}{\partial q} = q + w(t, X) \\
        \dot{q} = -\frac{\partial \boldsymbol{H}}{\partial x} = -(Dw(t, X))^\top q
    \end{cases}\label{eq:Hamsys}
\end{equation}
where $\boldsymbol{H}(x, q, t) = \frac{1}{2}|q|^2 +  q\cdot w(t, x)$ is the control Hamiltonian derived via Pontryagin’s minimum principle. Furthermore, if a potential $\varphi$ exists such that $q = \nabla \varphi$, then $\varphi$ satisfies the Hamilton-Jacobi equation:
\begin{equation*}
    \partial_t \varphi + \frac{1}{2}|\nabla \varphi|^2 + \nabla \varphi \cdot w = 0.
\end{equation*}
\end{theorem}

    The ODE system \eqref{eq:Hamsys} represents the bi-characteristic equations of the original PDE system \eqref{eq:optpdes}. The equation for $\dot{q}$ reveals that the optimal control velocity is not constant in the presence of a non-homogeneous medium; instead, it is transformed by the transpose of the Jacobian of the background flow, $(Dw)^T$. This ensures that the agent's path is a geodesic relative to the moving frame of the fluid.

\section{Optimization Strategy for the Initial Value Problem}\label{sec:optimization}

To solve for the optimal transport map in a moving medium, we adopt a shooting formulation at the level of individual agent trajectories. Let an agent's initial position $X(0) \in \mathbb{R}^d$ be drawn from the source distribution $\mu$, and we optimize over the control velocity $q_0$ to reach the target distribution $\nu$.

The minimizers of \eqref{eq:JHdist} subject to the boundary conditions $X(0)\sim \mu$ and $X(1)\sim \nu$ are realized by solving the following shooting problem:
\begin{align}
    \min_{q_0\in \mathbb{R}^d} \quad & \mathcal{F}\big( (X(1))_\# \mu , \nu \big) \label{eq:optotflowshooting}\\
    \text{s.t. } \quad & X(0) \sim \mu, q(0)=q_0,\nonumber \\
    & \dot X = q + w(t, X), \label{eq:odesysp_final} \\
    & \dot q = -(Dw(t, X))^\top q. \label{eq:odesysq_final}
\end{align}

The Hamiltonian system \eqref{eq:odesysp_final}--\eqref{eq:odesysq_final} provides the necessary conditions for optimality; any trajectory satisfying these ODEs is a stationary point of the control energy functional \eqref{eq:action}. The terminal cost functional $\mathcal{F}$ is utilized to enforce the boundary constraint $X(1) \sim \nu$. Common choices for $\mathcal{F}$ include the Kullback-Leibler (KL) divergence,
\begin{equation*}
    \mathcal{F}_{KL}(\rho, \nu) = \int_{\mathbb{R}^d} \rho(x) \log \left( \frac{\rho(x)}{\nu(x)} \right) dx,
\end{equation*}
or distance-based metrics such as the $L_p$ Wasserstein distance, the Maximum Mean Discrepancy (MMD) or the Sinkhorn divergence, which can be more robust in simulations for agent-based approximations.

\subsection{Linear Example}
Suppose the background flow is linear, $w(t, x) = Ax$, for a constant matrix $A \in \mathbb{R}^{d \times d}$. In this case, the Hamiltonian system \eqref{eq:odesysp_final}--\eqref{eq:odesysq_final} becomes:
\begin{align*}
  & \dot X = q + AX, \\
    & \dot q = -A^\top q
\end{align*}
For $\mu$-a.e. initial position $X(0)=\xi$, the control velocity has the closed-form solution $q(t) = e^{-A^\top t} q(0)$. Substituting this into the equation for $X(t)$ and using the variation of constants formula yields:
\begin{equation*}
    X(t) = e^{At} \left( \xi + \int_0^t e^{-As} e^{-A^\top s} q(0) \, ds \right) = e^{At} \left( \xi + \int_0^t e^{-(A + A^\top)s} q(0) \, ds \right).
\end{equation*}
At the terminal time $t=1$, the terminal value $X(1)$ satisfies:
\begin{equation}
    e^{-A} X(1) - \xi = \left( \int_0^1 e^{-(A + A^\top)s} \, ds \right) q(0). \label{eq:term_const_xv_lin}
\end{equation}
Defining the symmetric, positive-definite matrix $C = \int_0^1 e^{-(A + A^\top)s} \, ds$, we obtain the explicit shooting formula for the initial control velocity:
\begin{equation*}
    q(0) = C^{-1} \left( e^{-A} X(1) - \xi \right).
\end{equation*}

Since the cost \eqref{eq:K} is strictly convex in $q$, the dynamics are linear, and the terminal constraint \eqref{eq:term_const_xv_lin} is affine, the problem admits a unique global minimizer. The corresponding minimal energy is given by:
\begin{equation*}
    E_{\min} = \frac{1}{2} \int_{\mathbb{R}^d} (e^{-A} X(1) - \xi)^\top C^{-1} (e^{-A} X(1) - \xi) \, d\mu(\xi).
\end{equation*}
This expression corresponds to a static optimal transport problem between the measures $\mu$ and the pulled-back target $(e^{-A})_\# \nu$ with respect to the quadratic ground cost:
\begin{equation*}
    c_A(x, y) = \frac{1}{2}(e^{-A}y - x)^\top C^{-1} (e^{-A}y - x).
\end{equation*}
By Brenier's Theorem, there exists a Monge map $T$ such that the optimal assignment is given by $T(x) = x + C^{1/2} \nabla \psi(C^{-1/2} x)$ for some convex function $\psi$.

\subsection{Challenges in Optimization}
For a general background field $w(t, x)$, this formulation  \eqref{eq:optotflowshooting}$-$\eqref{eq:odesysq_final} is not equivalent to a static optimal transport problem. The presence of time-dependent background flows introduces two primary difficulties:
\begin{itemize}
    \item \textbf{Non-convexity:} The objective functional is generally non-convex with respect to the initial velocity $q_0$, and global minimizers are not guaranteed.
    \item \textbf{Coupled Matching:} The minimization is non-trivial because the point-to-point assignment between starting and target locations is not known a priori and must be discovered through the optimization of the initial velocity field. This is true even in the case of a linear background flow.
\end{itemize}

The numerical implementation of a shooting approach for solving \eqref{eq:optotflowshooting} is detailed in Algorithm \ref{alg:shooting_simple} in the next section.

\section{Numerical Implementation: Agent Discretization and Optimization} \label{sec:numericalimplementation}

While the continuum setting models the evolution of $X(t)$ starting from a single location, the numerical solution involves discretizing the measure $\mu$ into $N$ agents. This yields a system of $N$ coupled Hamiltonian ODEs for $\{X_i(t), q_i(t)\}_{i=1}^N$, where $X_i(t) \in \mathbb{R}^2$ denotes the position of agent $i$ and $q_i(t)$ represents its control velocity at time $t$.

We solve for the optimal initial controls $\mathbf{q}_0 = \{q_i(0)\}_{i=1}^N\in \mathbb{R}^{2N}$ by minimizing a regularized objective functional
\begin{equation}
    \mathcal{J}(\mathbf{q}_0) = \mathcal{F}_{\mathrm{KL}}(\hat{\rho}(1,\cdot), \nu) + \mathcal{P}_{\mathrm{boundary}}\label{eq:objective}
\end{equation}
using the L-BFGS quasi-Newton method \cite{Knoll2004Jacobian-freeApplications}. This functional consists of two terms described next: the first matches the empirical density of the agents, denoted by $\hat{\rho}$, to the target distribution $\nu$, and the second term $\mathcal{P}_{\mathrm{boundary}}$ ensures that agents stay within the computational domain.

\subsection{Agent Dynamics and Density Estimation}
The agents evolve according to the Hamiltonian system \eqref{eq:odesysp_final}--\eqref{eq:odesysq_final}, integrated using an adaptive explicit Runge--Kutta method (RK45) \cite{ShampineSomeFormulas, DormandAFormulae}. The initial positions are sampled such that $\hat{\rho}(0, \cdot) \sim \mu$, and the initial velocities are assigned from the optimization variable $\mathbf{q}_0$.

To evaluate the KL divergence, we approximate the evolving density field $\rho(t, x)$ using the empirical measure:
\begin{equation*}
    \hat{\rho}(t, x) = \frac{1}{N} \sum_{i=1}^N \mathcal{K}_\sigma(x - X_i(t)),
\end{equation*}
where $\mathcal{K}_\sigma$ is a kernel density estimator (KDE)  with bandwidth $\sigma$ centered at the agent locations \cite{Silverman2018DensityAnalysis}. This regularization transforms the sum of Dirac deltas into a smooth, differentiable density function suitable for gradient-based optimization.

\subsection{Boundary Regularization}
To ensure the trajectories remain within a physically meaningful domain, we introduce a penalty term $\mathcal{P}_{\mathrm{boundary}}$ for agents whose coordinates exceed a prescribed radius $D > 0$. We define the violation function:
\begin{equation*}
    \phi(X_i(t)) = \max(0, \|X_i(t)\| - D).
\end{equation*}
The boundary penalty is computed both at the terminal time and as an aggregate over the entire trajectory to prevent transient excursions:
\begin{equation*}
    \mathcal{P}_{\mathrm{boundary}} = \lambda_b \sum_{i=1}^N \phi(X_i(1))^2 + \frac{\lambda_b}{T} \int_0^1 \sum_{i=1}^N \phi(X_i(t))^2\,dt,
\end{equation*}
where $\lambda_b$ is a penalty parameter. This formulation ensures that the shooting method favors trajectories that stay within the domain while seeking to match the target distribution $\nu$.

\begin{figure}[t]
\centering
\captionsetup[subfigure]{labelformat=simple} 
\renewcommand\thesubfigure{(\alph{subfigure})} 
    \begin{subfigure}[b]{0.32\textwidth}
        \centering
        \includegraphics[width=\textwidth]{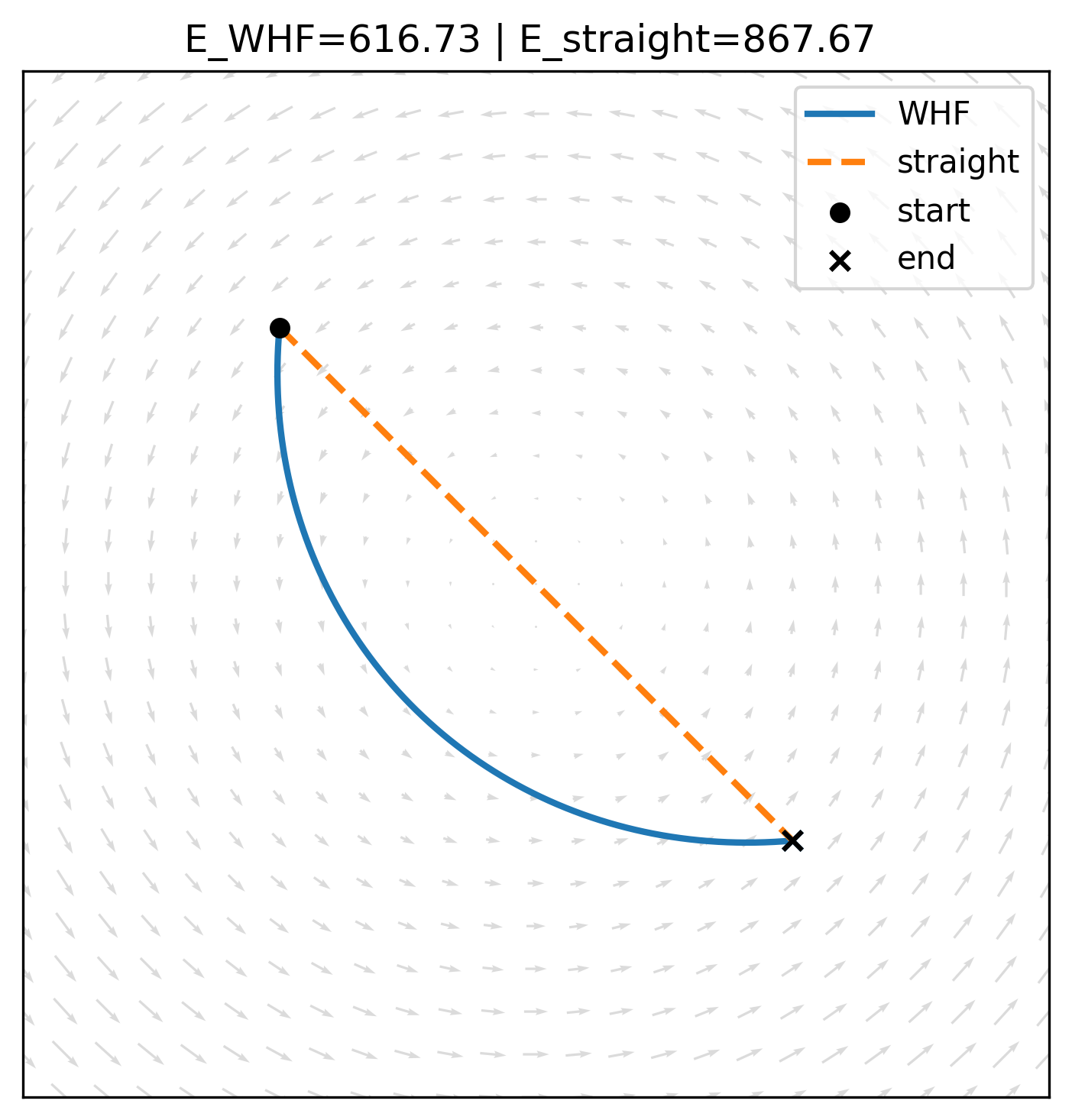}
        \caption{Circle}
    \end{subfigure}
    \hfill
    \begin{subfigure}[b]{0.32\textwidth}
        \centering
        \includegraphics[width=\textwidth]{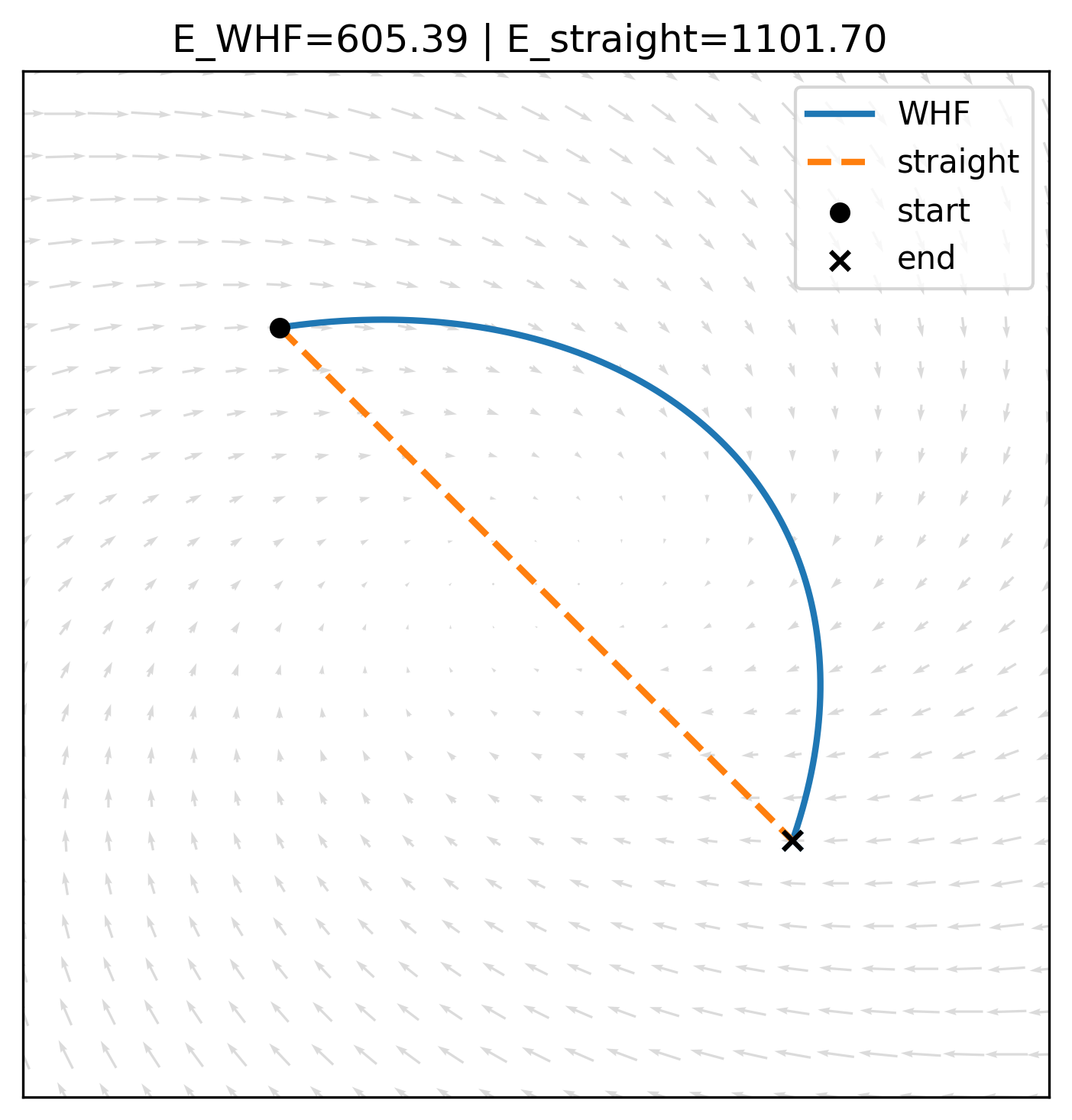}
        \caption{Attractor}
    \end{subfigure}
    \hfill
    \begin{subfigure}[b]{0.32\textwidth}
        \centering
        \includegraphics[width=\textwidth]{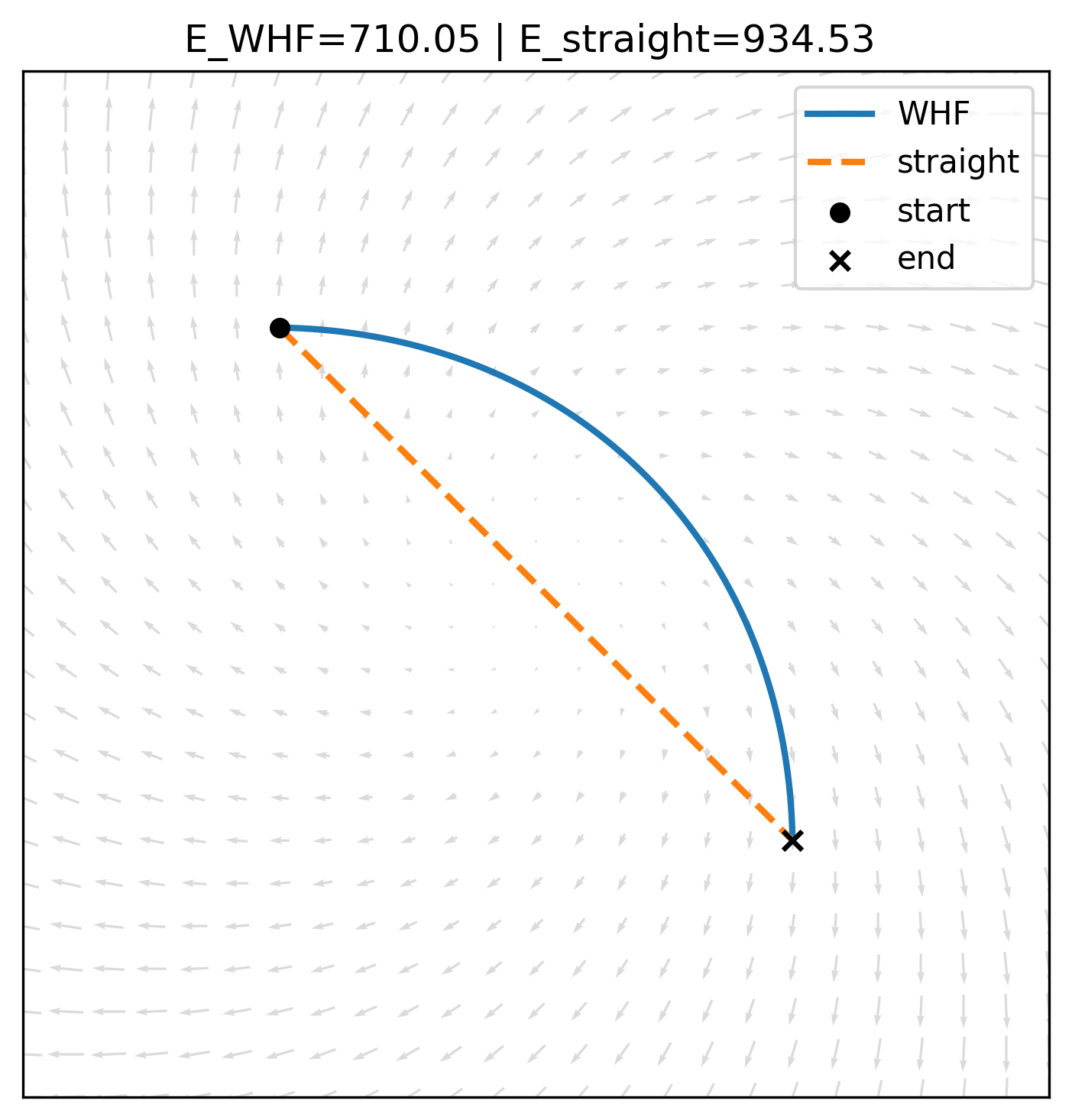}
        \caption{Repeller}
    \end{subfigure}

    \vspace{0.5cm} 

    \begin{subfigure}[b]{0.32\textwidth}
        \centering
        \includegraphics[width=\textwidth]{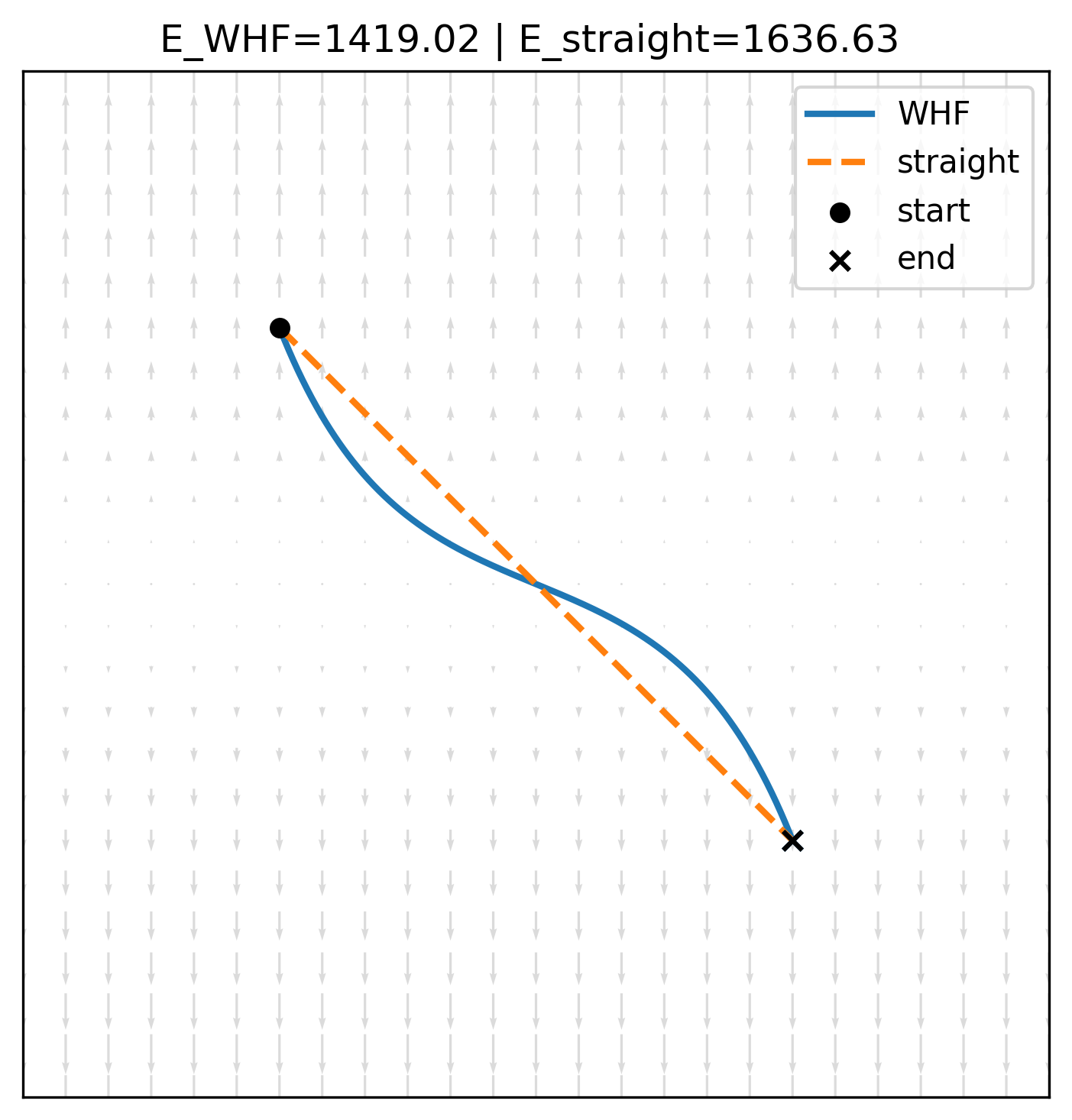}
        \caption{Vertical}
    \end{subfigure}
    \hfill
    \begin{subfigure}[b]{0.32\textwidth}
        \centering
        \includegraphics[width=\textwidth]{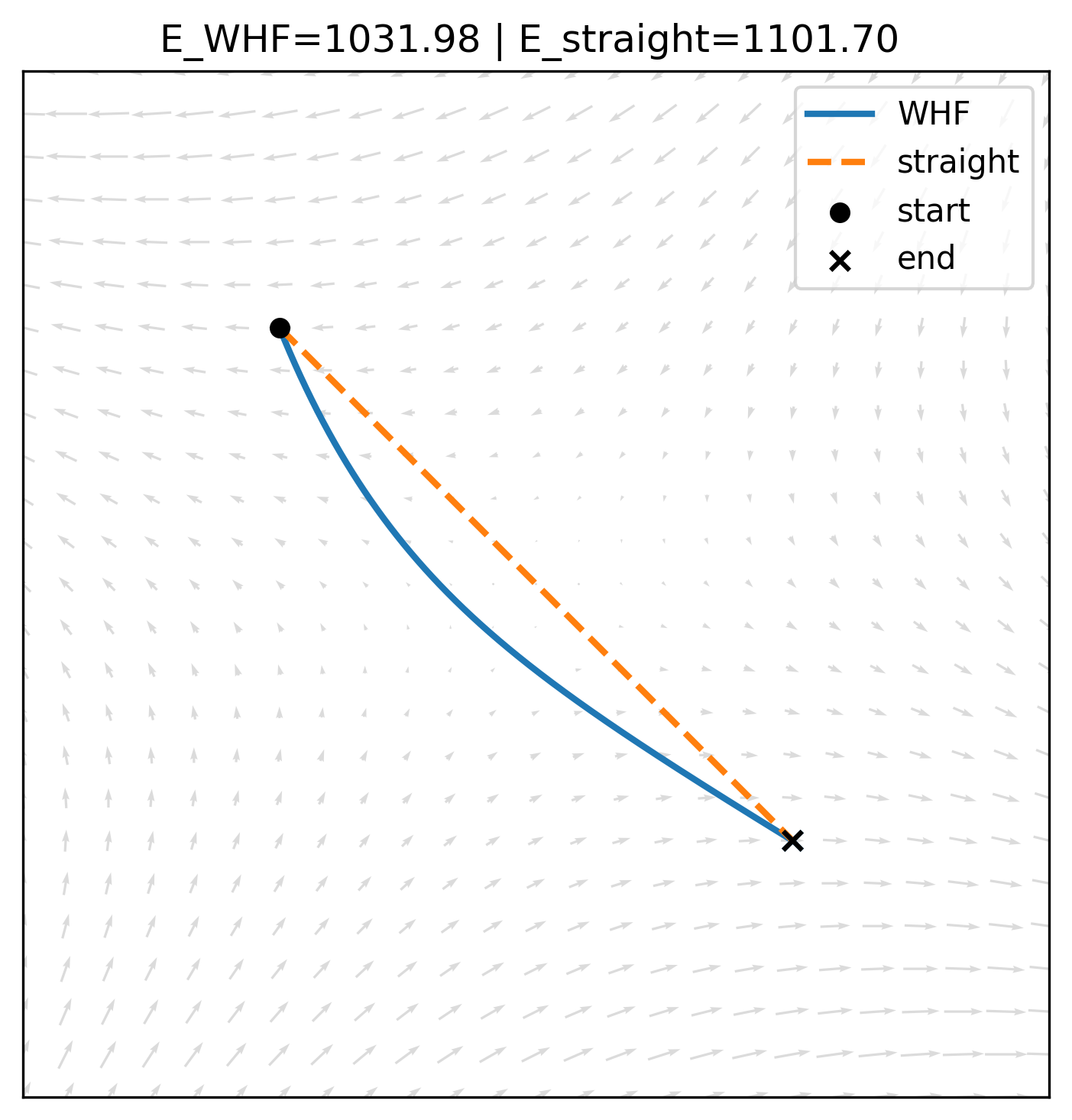}
        \caption{Stagnation}
    \end{subfigure}
    \hfill
    \begin{subfigure}[b]{0.32\textwidth}
        \centering
        \includegraphics[width=\textwidth]{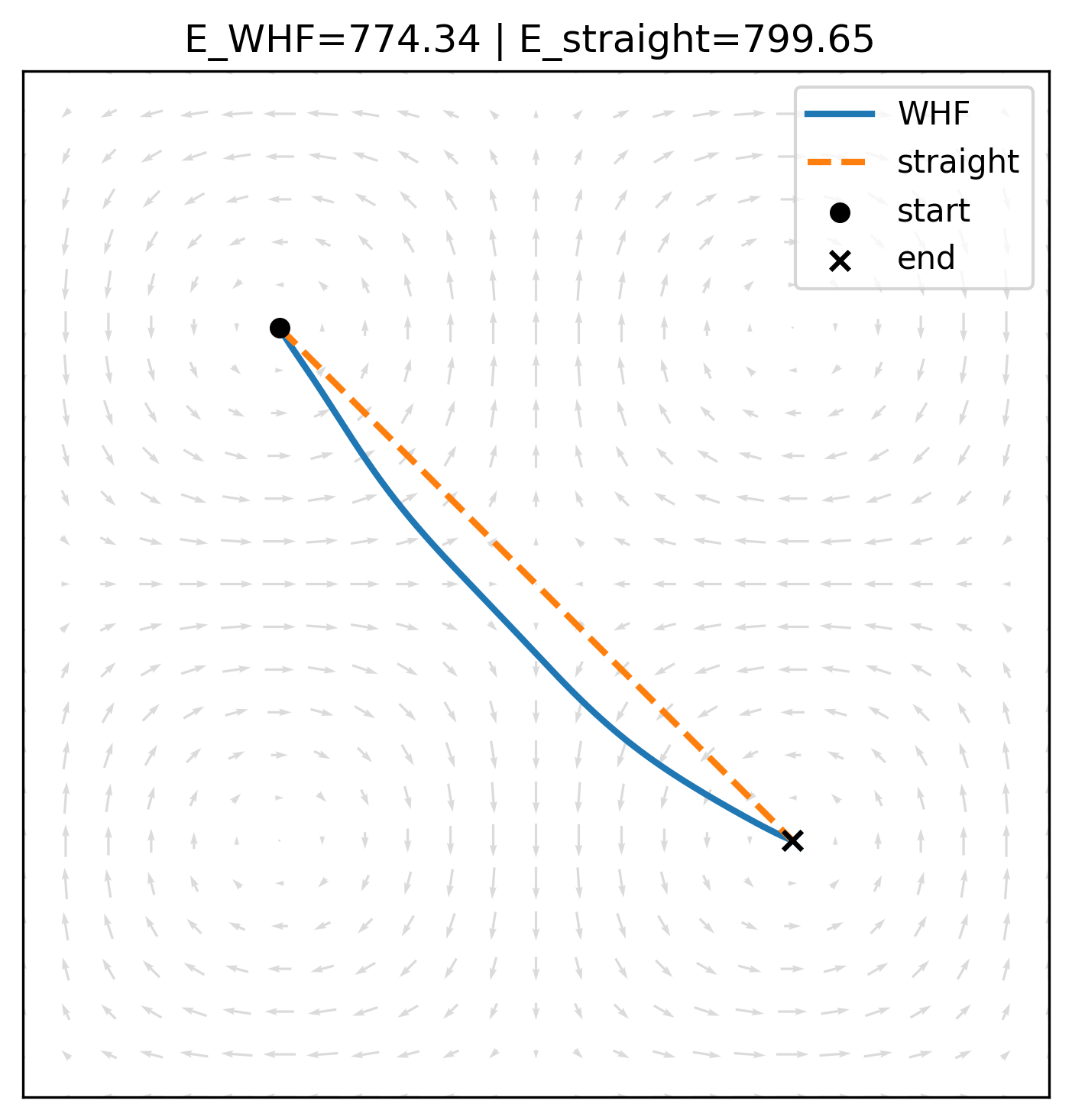}
        \caption{Gyre}
    \end{subfigure}

\caption{Single-agent path planning in representative flow fields. We consider several steady and time-dependent background flows:
(a) Circle: $w(x,y) = (-y,\,x)^\top$, 
(b) Point Attractor: $w(x,y) = (-x+2y,\,-y-x)^\top$,
(c) Point Repeller: $w(x,y) = (x+y,\,-x+y)^\top$,
(d) Vertical: $w(x,y) = (0,\,5y)^\top$,
(e) Stagnation Point: $w(x,y) = (x-2y,\,-x-y)^\top$, and
(f) Gyre: $w(t,x,y) =\left(-2\pi\sin(\pi f),\, 2\pi\cos(\pi f)\,\partial_x f \right)^\top$, 
where \( f(t,x) = a(t)x^2 + b(t)x \), \( a(t) = \epsilon\sin(\omega t) \), and \( b(t) = 1 - 2\epsilon\sin(\omega t) \), with \( \epsilon = 0.1 \), \( \omega = 2\pi \).}
\label{fig:otflows1N}
\end{figure}

\begin{algorithm}
\caption{Shooting-based Optimal Transport in a Moving Medium}
\label{alg:shooting_simple}
\begin{algorithmic}[1]
\Require Source density $\mu$, target density $\nu$, background flow field $w(t, x)$
\Require Time horizon $T=1$, agent count $N$, KDE bandwidth $\sigma$, penalty weight $\lambda_b$
\Ensure Optimal initial velocities $\mathbf{q}_0$, agent trajectories $\mathbf{X}(t)$

\State Sample initial agent positions $\{X_i(0)\}_{i=1}^N$ from $\mu$
\State Initialize initial control velocities $\mathbf{q}_0 \in \mathbb{R}^{N \times d}$ \Comment{e.g., zero initialization or warm-start}

\Function{Objective}{$\mathbf{q}_0$} \Comment{Evaluates $\mathcal{J}(\mathbf{q}_0)$ defined in \eqref{eq:objective}}
    \State Integrate the Hamiltonian system for $t \in [0, 1]$ using RK45:
    \[
    \begin{cases} 
    \dot{X}_i = q_i + w(t, X_i) \\
    \dot{q}_i = -(Dw(t, X_i))^\top q_i 
    \end{cases}
    \]
    \State Compute the empirical density $\hat{\rho}(1, \cdot)$ at terminal time using KDE with bandwidth $\sigma$
    \State $\mathcal{L}_{\text{match}} \gets \mathcal{F}_{\mathrm{KL}}(\hat{\rho}(1, \cdot), \nu)$ \Comment{Target matching loss}
    \State $\mathcal{P}_{\text{bound}} \gets \text{Compute boundary penalty via violation function } \phi$
    \State \Return $\mathcal{J} = \mathcal{L}_{\text{match}} + \mathcal{P}_{\text{bound}}$
\EndFunction

\State Minimize $\mathcal{J}(\mathbf{q}_0)$ with respect to $\mathbf{q}_0$ using L-BFGS \Comment{Quasi-Newton optimization}

\State \Return $\mathbf{q}_0^*$ and the resulting trajectories $\{X_i(t)\}_{i=1}^N$
\end{algorithmic}
\end{algorithm}

\section{Numerical Experiments}\label{sec:numerics}

Equations \eqref{eq:odesysp_final}–\eqref{eq:odesysq_final} are integrated over $[0, T=1]$ using an adaptive explicit Runge--Kutta method (RK45) with absolute tolerance \( 10^{-9} \) and timestep $\Delta t=.001$ unless otherwise stated. The density estimation and KL divergence are evaluated on a uniform spatial grid of size $500\times 500$ over the domain \([-20, 20]^2\). The numerical experiments were conducted on a laptop featuring an Apple M3 Pro processor (12-core CPU) and 18 GB of unified memory running macOS 14.4.1. All algorithms were implemented in Python 3.14. Linear algebra operations were performed using NumPy, linked against the Apple Accelerate framework for optimized BLAS and LAPACK routines. The code was compiled using Clang 15.0.0 and utilized Advanced SIMD (ASIMD) extensions for vectorized computations.

Optimization is carried out via L--BFGS (SciPy), with convergence tolerance \(10^{-6}\). All simulations are performed for ensemble sizes \( N \in \{1, 5, 10, 25, 50\} \). For \( N = 1 \), the setup corresponds to single-agent planning (Figure~\ref{fig:otflows1N}).

The effort is estimated from the $N$ trajectories $X_i$ via
\[
E(N) = \sum_{i=1}^N\int_0^T \|\dot{X_i}(t) - w(t, X_i(t))\|^2\,dt,
\]
using central finite differences to compute $\dot{X}_i$. Energy savings are computed relative to straight-line baselines:
\[
\mathrm{Savings}(N) = 1 - \frac{E_{\mathrm{WHF}}(N)}{E_{\mathrm{str}}(N)}.
\]

\subsection{Single-agent experiments (Gaussian target)}

We choose a target density consisting of a single isotropic Gaussian centered at 
$\boldsymbol{c}=(c_x,c_y)$ with width~$s$
\[\displaystyle 
\nu_{\textrm{point}}(x,y) \propto \frac{1}{2\pi s^2}\exp\!\Big(-\frac{\|(x,y) - \boldsymbol{c}\|^2}{2s^2}\Big).\]

We first visualize the transport of a single agent under each background flow (Figure \ref{fig:otflows1N}).
The initial position is \(X(0) = (-10,10)^\top\) and the target density center is \({\bf c}=(10,-10)^\top\) and the KDE bandwidth for the agent density is $\sigma=1$ and $s=10$.
Figure~\ref{fig:otflows1N} shows the optimized path (solid blue) overlaid with the straight-line trajectory (dashed orange) and the flow field (gray quivers).
For reference, the total efforts \(E_{\mathrm{opt}}\) and \(E_{\mathrm{str}}\) are reported in Table~\ref{tab:single_agent_flows}.
Flows such as the point attractor, point repeller, and circle substantially reduce the required effort due to the alignment of the flow with the direction of travel. In the other cases, the agent must spend more effort ``fighting against the current,'' and therefore the savings are not as dramatic in these cases. The optimization routines took less than 7 seconds to complete within 2-13 iterations.

\begin{table}[ht]
\centering
\begin{tabular}{lrrrrr}
\hline
flow type & $E_{\mathrm{WHF}}$ & $E_{\mathrm{str}}$ & Savings(\%) & runtime (s) & iterations \\
\hline
Circle          & 616.73  & 867.67  & 28.92 & 0.99 & 2  \\
Attractor     & 605.39  & 1101.70 & 45.05 & 5.53 & 11 \\
Repeller   & 710.05  & 934.53  & 24.02 & 1.01 & 2  \\
Vertical        & 1419.02 & 1636.64 & 13.30 & 6.23 & 13 \\
Stagnation     & 1031.98 & 1101.70 & 6.33  & 5.23 & 12 \\
Gyre            & 774.34  & 799.65  & 3.17  & 6.01 & 6  \\
\hline
\end{tabular}
\caption{Performance metrics for single-agent path planning across various background flows ($N=1$). $E_{\mathrm{WHF}}$ denotes the control effort (energy) using the WHF method, while $E_{\mathrm{str}}$ represents the energy required for a straight-line trajectory baseline. Savings indicate the percentage reduction in control effort achieved by exploiting flow geometry. Runtimes and iteration counts reflect the efficiency of the L-BFGS-B optimizer.}
\label{tab:single_agent_flows}
\end{table}

\begin{figure}[h!]
    \centering
    \includegraphics[width=\linewidth]{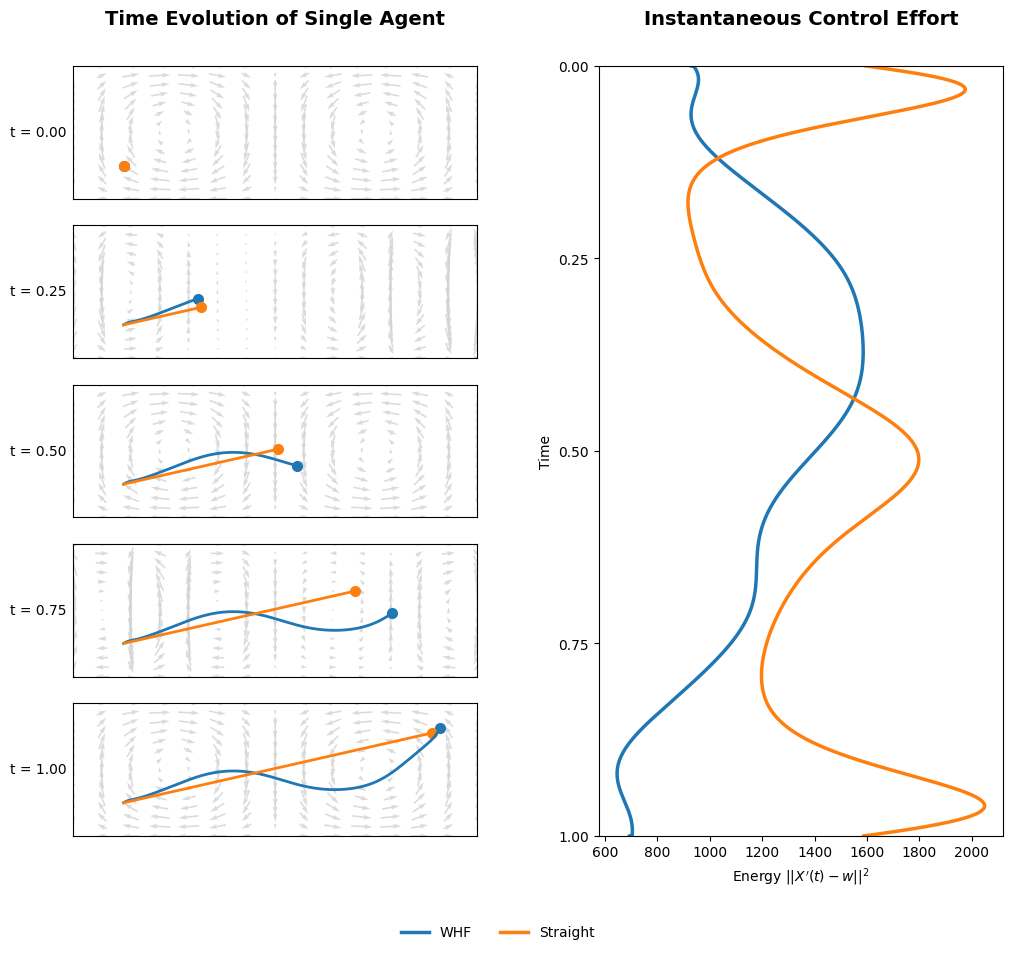}
    \caption{Left: Single-agent path-planning in a time-dependent `gyre' background flow: $w(t,x,y) =\left(-2\pi\sin(\pi f),\, 2\pi\cos(\pi f)\,\partial_x f \right)^\top$, 
where \( f(t,x) = a(t)x^2 + b(t)x \), \( a(t) = \epsilon\sin(\omega t) \), and \( b(t) = 1 - 2\epsilon\sin(\omega t) \), with \( \epsilon = 0.1 \), \( \omega = 2\pi \). Right: Instantaneous control effort $\|\dot{X}(t) - w(t, X(t))\|^2,$. }
    \label{fig:gyre}
\end{figure}

Figure \ref{fig:gyre} shows five time snapshots of the single-agent path-planning algorithm (Algorithm \ref{alg:shooting_simple}) in a time-dependent flow, and a plot of the instantaneous control effort during the environment traversal. The snapshots demonstrate how the optimal path aligns with the background flow at each time step, resulting in energy savings.
The right plot in Figure \ref{fig:gyre} shows that the instantaneous control energy $\|\dot{X}(t) - w(t, X(t))\|^2,$ is not constant over time (unlike in classical optimal transport). The optimal control energy is higher when the agents move against the current and lower when the best directions align with the flow field.

\subsection{Multi-agent ensembles}

To evaluate the robustness of the  WHF method against varying agent counts and target geometries, we perform a series of numerical sweeps. The simulation environment is defined on a two-dimensional domain $\Omega \in [-20, 20]^2$, discretized by a $500 \times 500$ grid. The temporal domain is defined by a horizon $T = 1.0$ and a step size $\Delta t = 0.001$ ($1000$ integration steps).

For each case, the agents are initially positioned in a circular formation of unit radius centered at the origin:
\begin{equation}
    \mathbf{p}_{0,i} = \begin{bmatrix} \cos \theta_i \\ \sin \theta_i \end{bmatrix}, \quad \theta_i = \frac{2\pi(i-1)}{N}.
\end{equation}

The target density is an annular target of mean radius~$r_0$ centered at~$\boldsymbol{c}=(0,0)^\top$ and 
width controlled by~$s$: 

\[\displaystyle
\nu_{\textrm{ring}}(x,y)  \propto
\exp\!\Big(-\frac{(R - r_0)^2}{2s^2}\Big), 
\quad R=\sqrt{(x-c_1)^2+(y-c_2)^2}\]

\begin{figure}[ht]
     \centering
     \begin{subfigure}[b]{0.5\textwidth}
         \centering
         \includegraphics[width=\linewidth,trim={0 0 3.5cm 0}, clip]{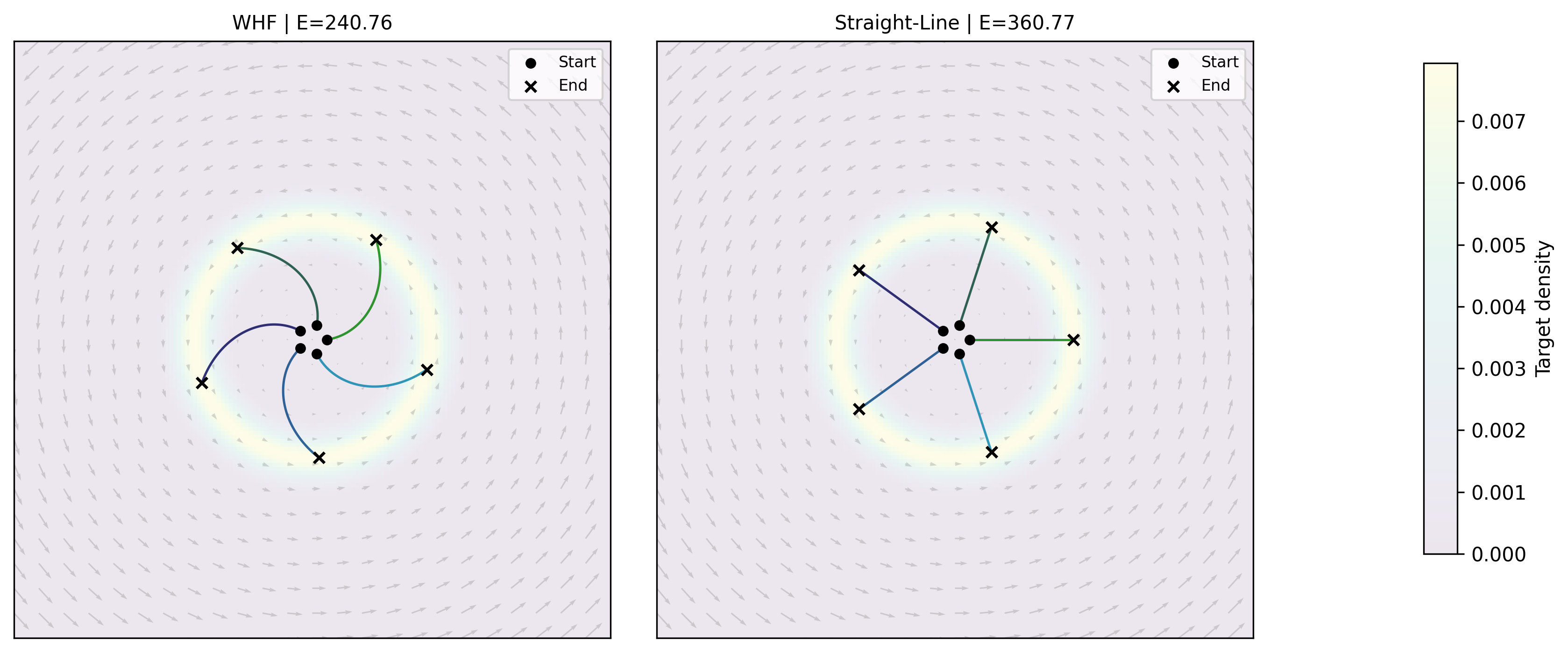}
         \caption{$N=5$}
         \label{fig:n5}
     \end{subfigure}\begin{subfigure}[b]{0.5\textwidth}
         \centering
         \includegraphics[width=\linewidth,trim={0 0 3.5cm 0}, clip]{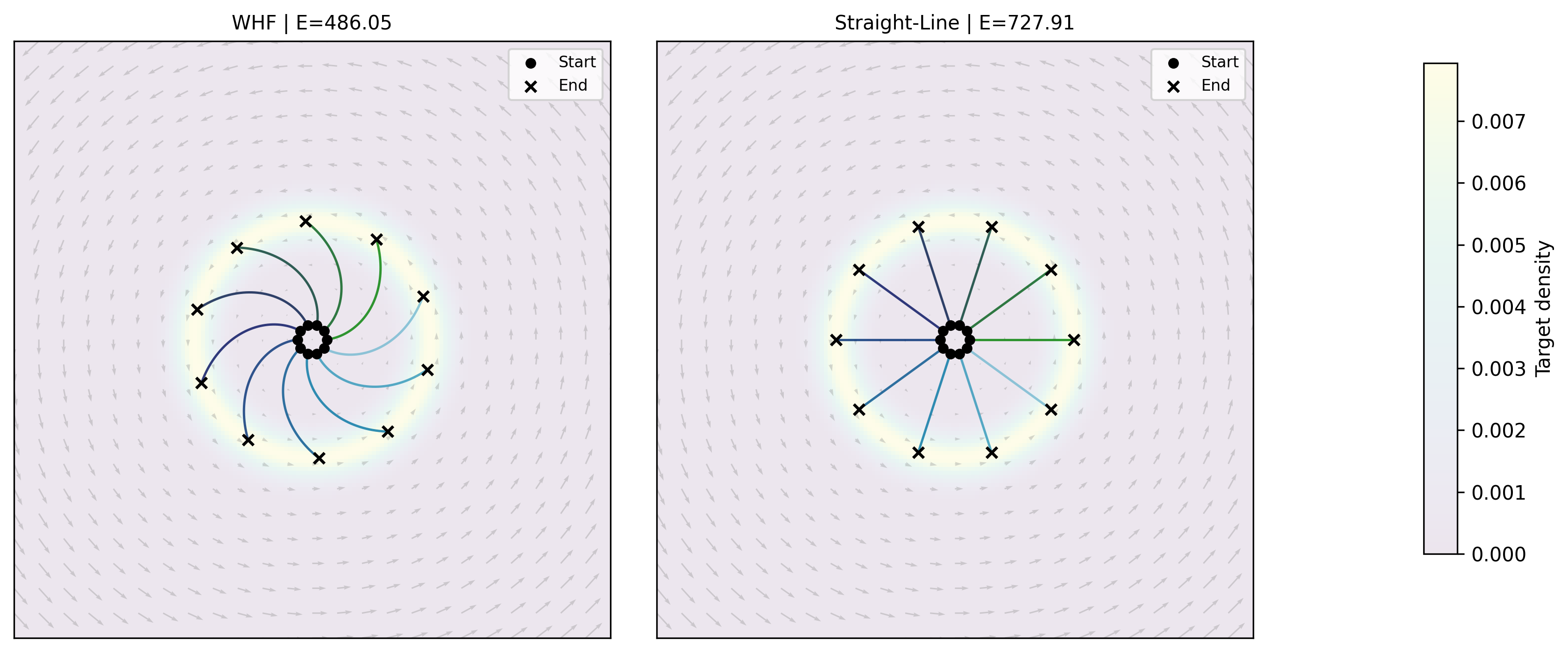}
         \caption{$N=10$}
         \label{fig:n10}
     \end{subfigure}
     
     \vspace{0.25cm} 

     \begin{subfigure}[b]{0.5\textwidth}
         \centering
         \includegraphics[width=\linewidth,trim={0 0 3.5cm 0}, clip]{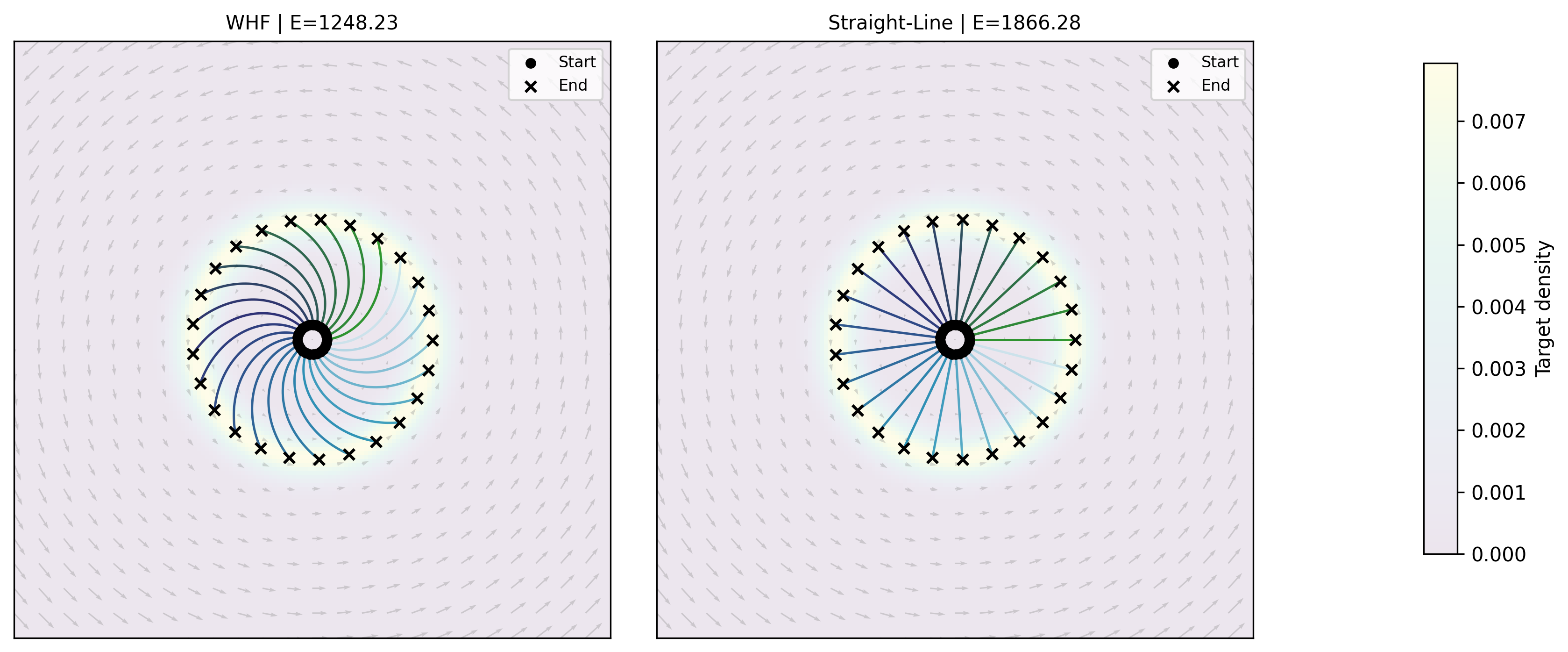}
         \caption{$N=25$}
         \label{fig:n25}
     \end{subfigure}\begin{subfigure}[b]{0.5\textwidth}
         \centering
         \includegraphics[width=\linewidth,trim={0 0 3.5cm 0}, clip]{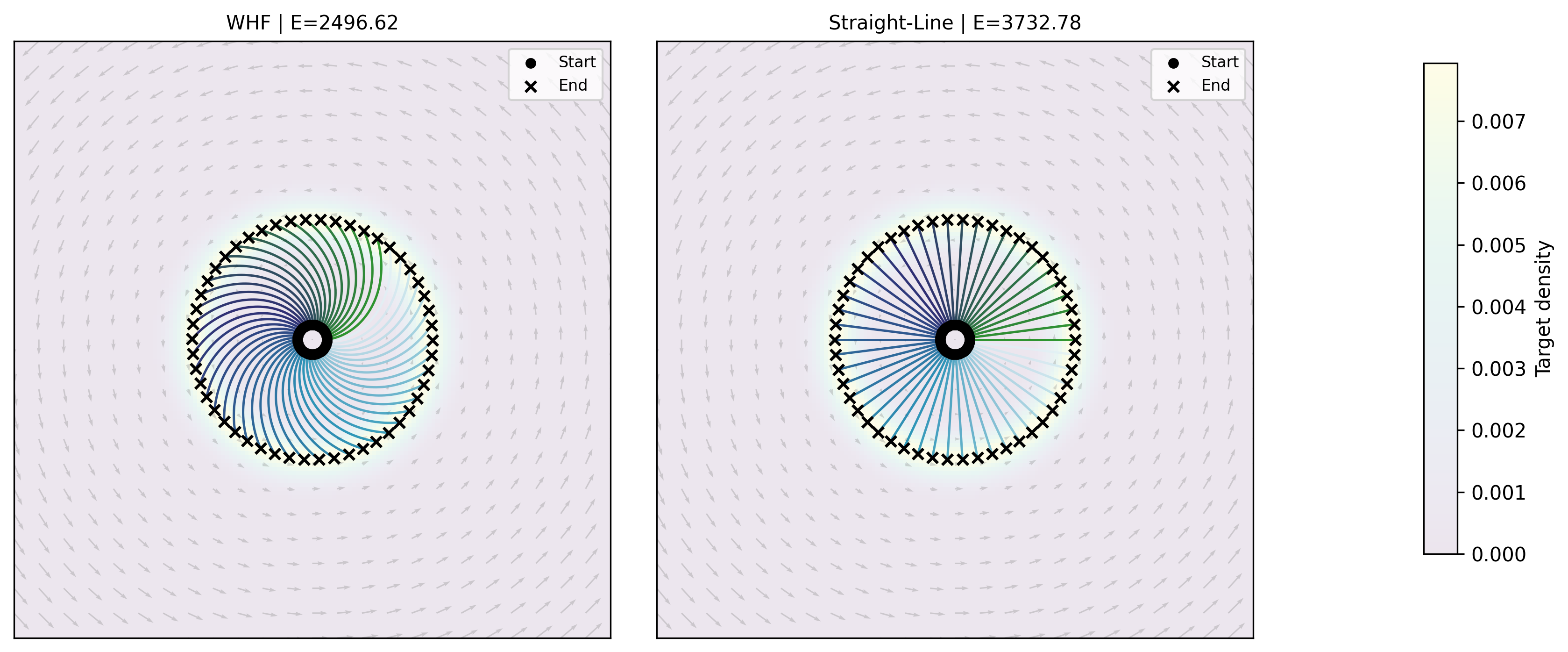}
         \caption{$N=50$}
         \label{fig:n50}
     \end{subfigure}
     
    \caption{Multi-agent trajectory optimization and target formation in a rotational ``circle'' background flow for varying ensemble sizes. The agents navigate from a compact initial cluster at the origin to an equilateral distribution along the annular target manifold. The results demonstrate that the Wasserstein-Hamiltonian Flow (WHF) maintains high-fidelity density matching via KL-divergence minimization while scaling effectively with agent count $N$.}
    \label{fig:multiagent}
\end{figure}

Plot \ref{fig:multiagent} shows the multiagent target formation in the case of a `circle' background flow. In each case, the final locations are evenly spread out on the target ring due to the minimization of the KL-divergence between the final distribution and target distribution.


\begin{figure}

\includegraphics[width=.8\linewidth]{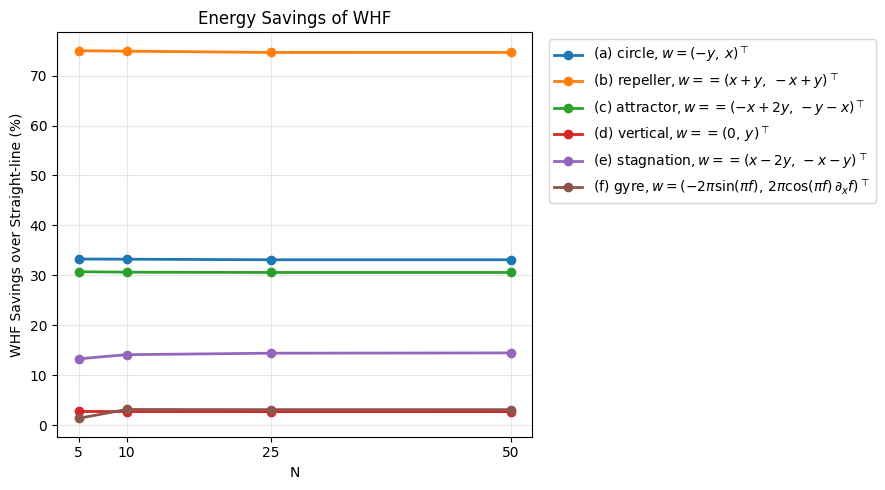}
\caption{The savings in energy are proportional to the number of agents, as demonstrated by these plots.}
\end{figure}



\subsubsection{Local minima}

Complex target geometries often introduce non-convexities in the loss landscape, causing the optimization to converge to sub-optimal local minima. Here, we consider a background \textit{point repeller} flow field. The temporal domain is defined by a horizon $T = 1$ and a step size $\Delta t = 0.01$ ($100$ integration steps). We analyze two agent populations, $N \in \{10, 25\}$. Two distinct target distributions are evaluated: a symmetric \textit{ring} distribution ($s=1, r_0=8$) and a non-symmetric \textit{heart} distribution, defined as:

\begin{equation}
    \nu_{\text{heart}}(x, y) \propto \exp\left( -\frac{\Phi_{\text{heart}}(x, y)}{2s^2} \right),
    \label{eq:heart_density}
\end{equation}

\noindent where the shape-defining potential $\Phi_{\text{heart}}$ is given by:

\begin{equation}
    \Phi_{\text{heart}}(x, y) = X_l^2 + \left( \frac{5}{4}Y_l - \sqrt{|X_l|} \right)^2.
\end{equation}

\noindent The coordinates $(X_l, Y_l)$ represent the centered and scaled spatial domain:

\begin{equation}
    \begin{aligned}
        X_l &= l(x - c_x), \\
        Y_l &= l(y - c_y),
    \end{aligned}
\end{equation}

\noindent where $\mathbf{c} = (c_x, c_y)$ is the target centroid and $u$ is a reciprocal scaling factor that determines the spatial extent of the distribution. In this study, we set $s=3$ to control the boundary thickness and $l=0.15$ to scale the target within the computational domain $\Omega = [-20, 20]^2$.

To assess the convergence of the optimizer in the presence of the background flow, we conduct $M=100$ trials for each $(N, \nu_{\text{target}})$ pair. The initial velocity $\mathbf{q}_0$ for the optimization is drawn from a uniform distribution:
\begin{equation}
    \mathbf{q}_{0,i} \sim \mathcal{U}(-1, 1)^2.
\end{equation}
For each trial, we record the total energy cost $E_{\mathrm{WHF}}$ and the final value of the cost function $\mathcal{J}(\mathbf{q}_{0,\mathrm{opt}})$ upon reaching the stopping criteria (maximum $300$ iterations).

\begin{figure}
    \centering
        \includegraphics[width=0.3\linewidth]{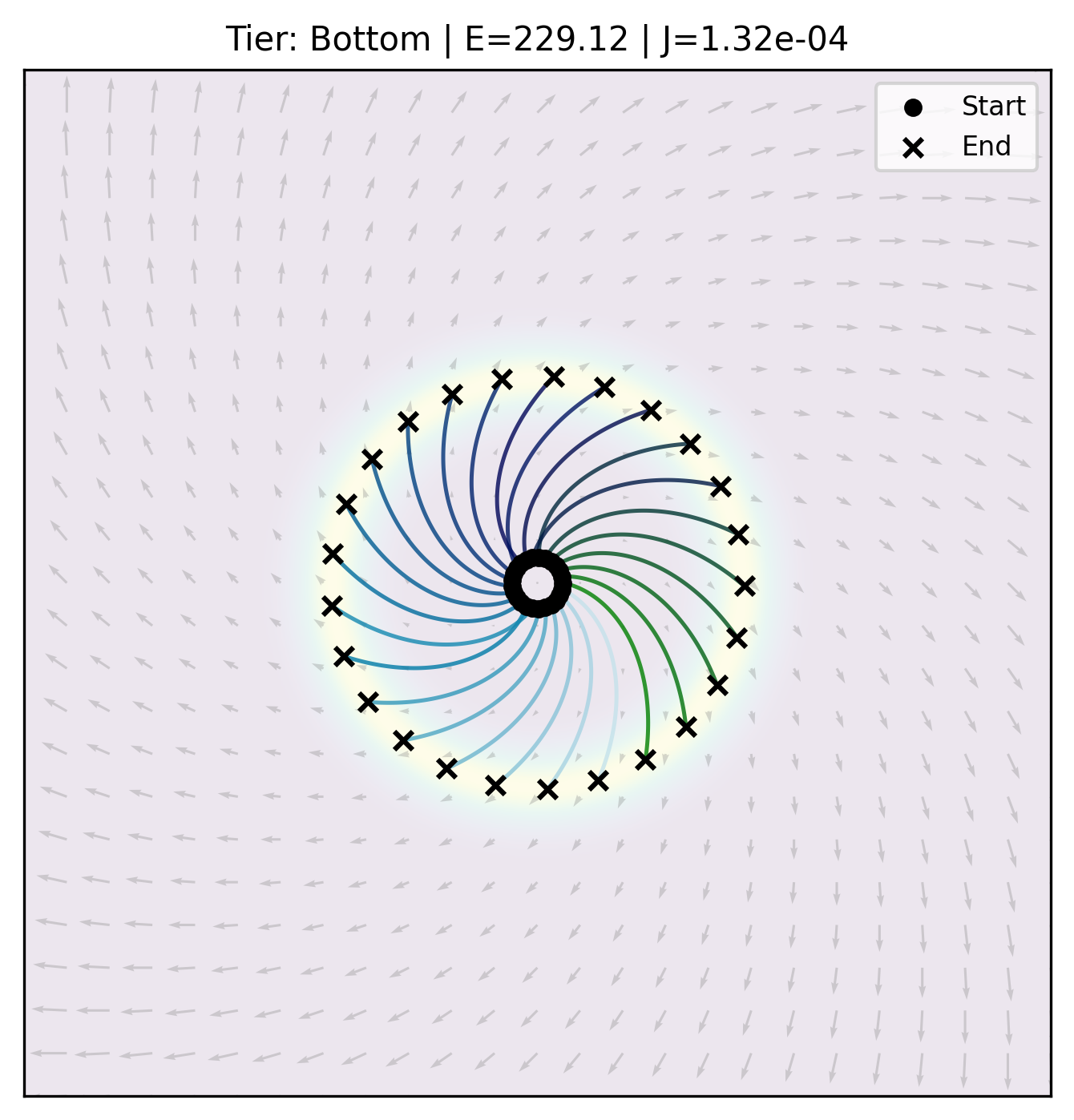}  \includegraphics[width=0.3\linewidth]{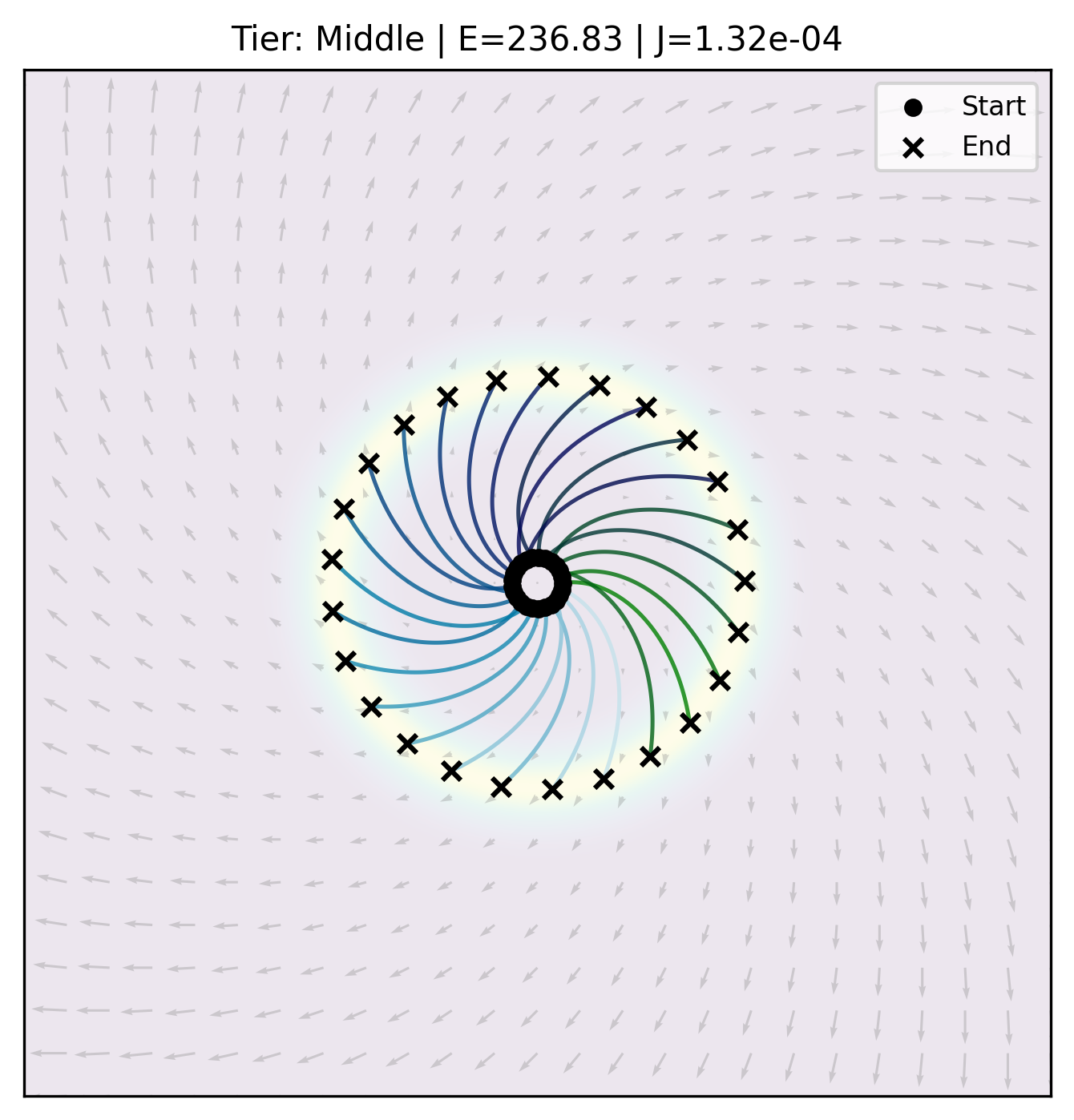}  \includegraphics[width=0.3\linewidth]{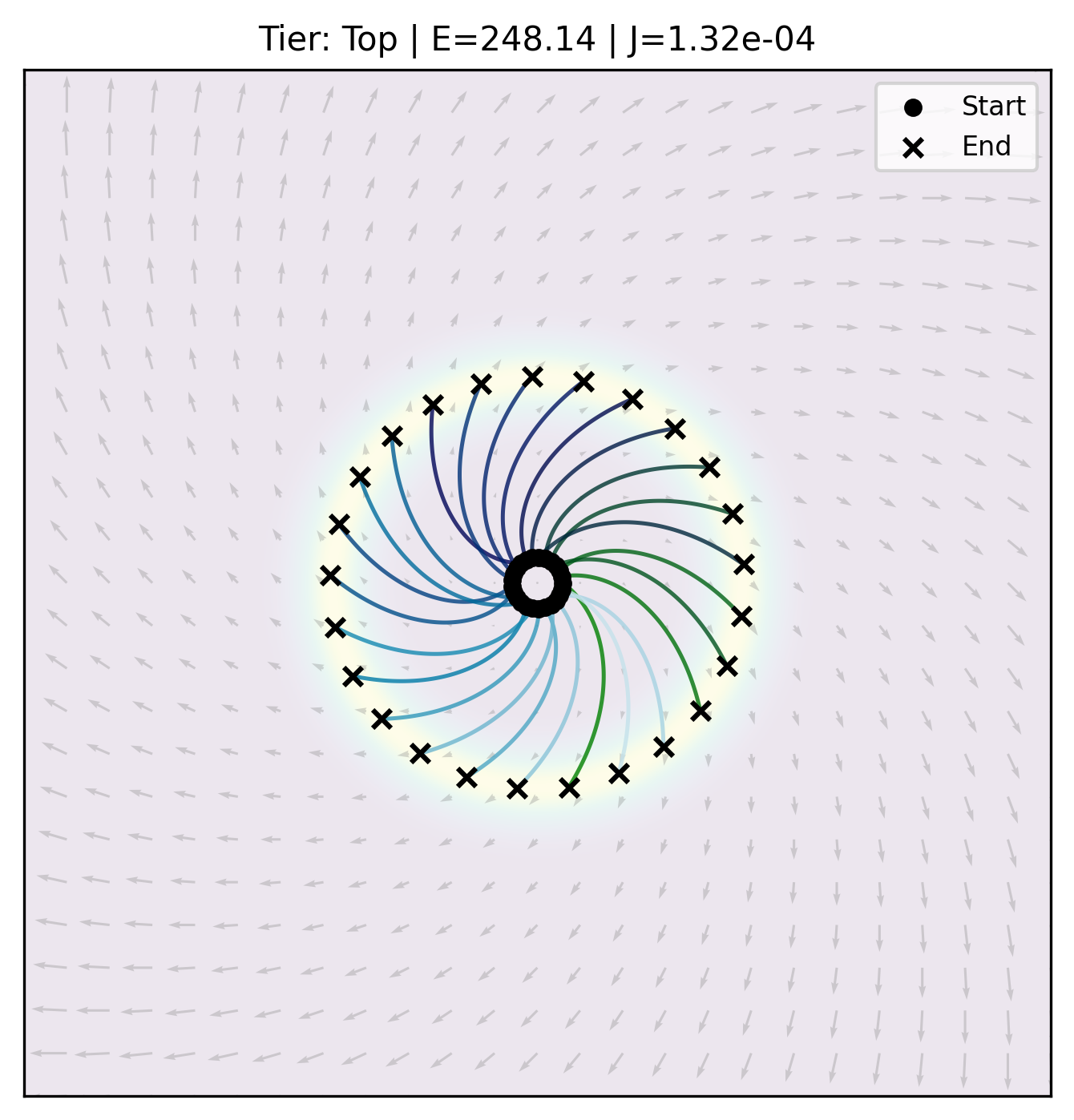}
     \includegraphics[width=0.3\linewidth]{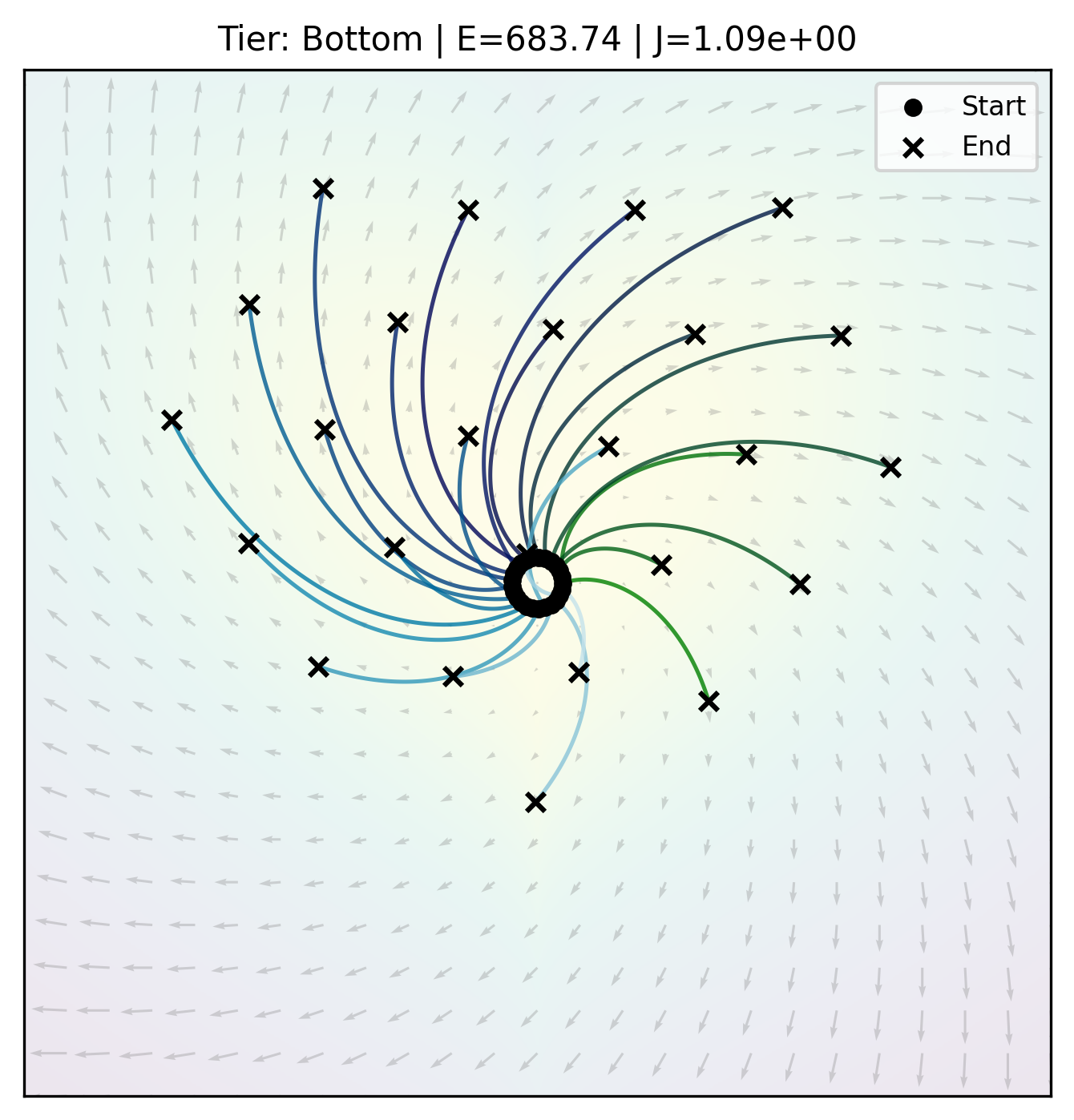}  \includegraphics[width=0.3\linewidth]{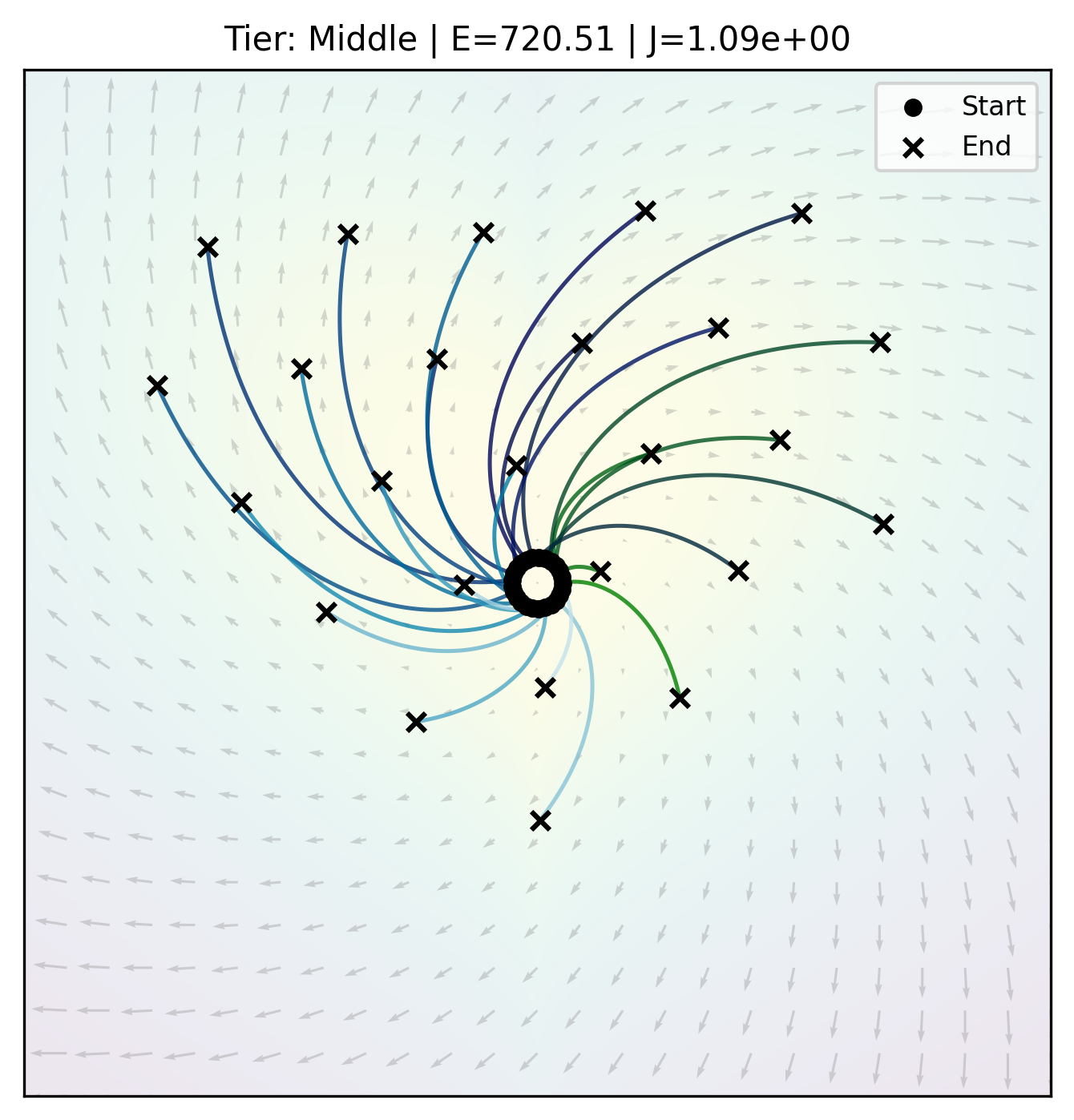}  \includegraphics[width=0.3\linewidth]{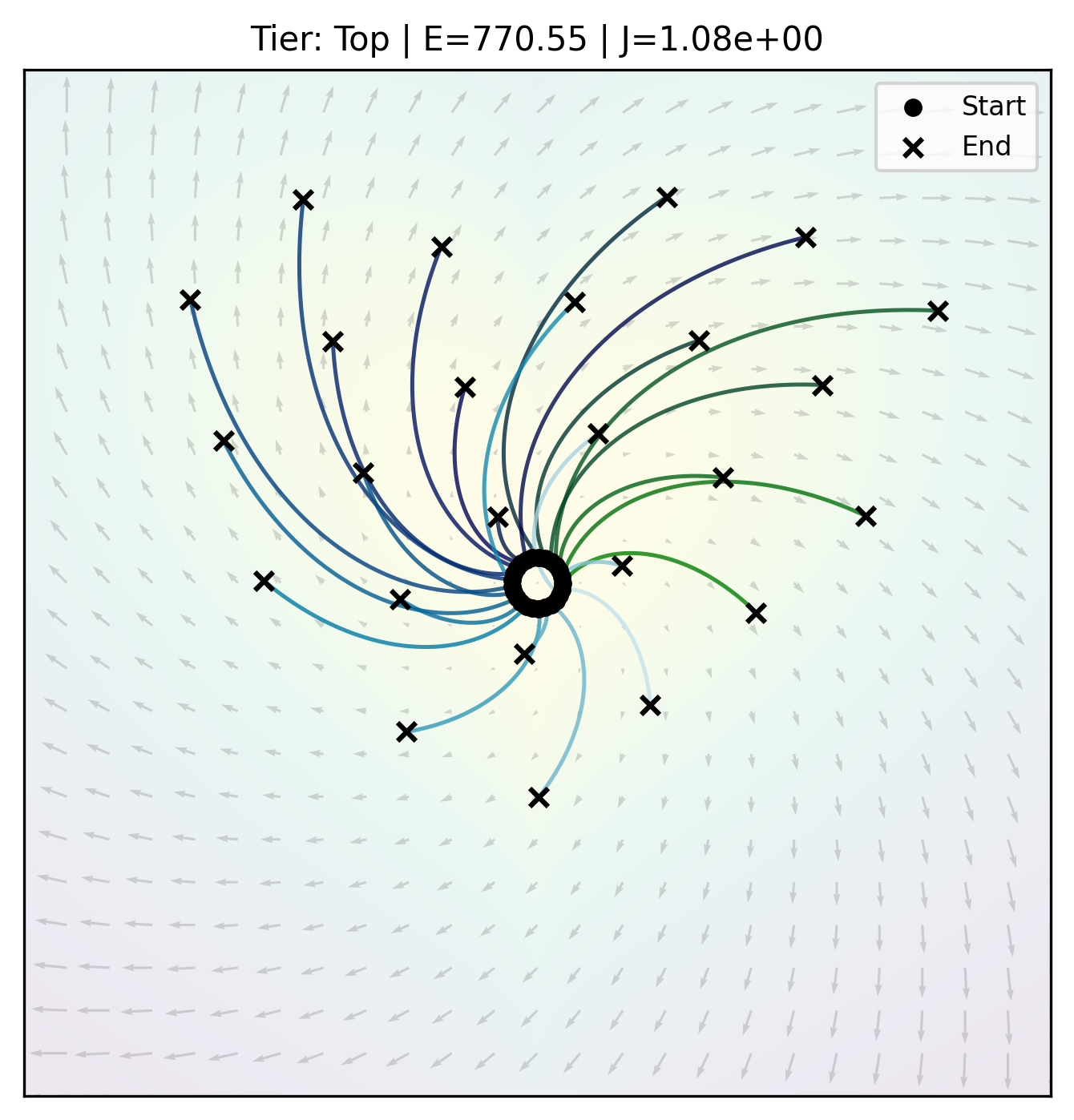}
    \caption{Plots generated by minimizers of \eqref{eq:objective} for the target distributions $\nu_{\mathrm{ring}}$ (top) and $\nu_{\mathrm{heart}}$ (bottom) and $N=25$ agents. The final objective values are similar, however they have different control energy values, each corresponding to local minimizers. The plots show energy expenditure in increasing order from left to right.}
    \label{fig:localmin}
\end{figure}

Figure \ref{fig:localmin} shows plots generated by minimizers of \eqref{eq:objective} for the target distributions $\nu_{\mathrm{ring}}$ (top) and $\nu_{\mathrm{heart}}$ (bottom) and $N=25$ agents. The final objective values are similar, however they have different control energy values, each corresponding to local minimizers. The plots show energy expenditure in increasing order from left to right. 

\begin{figure}[t]
    \centering
    \begin{subfigure}{0.48\textwidth}
        \centering
        \includegraphics[width=\linewidth]{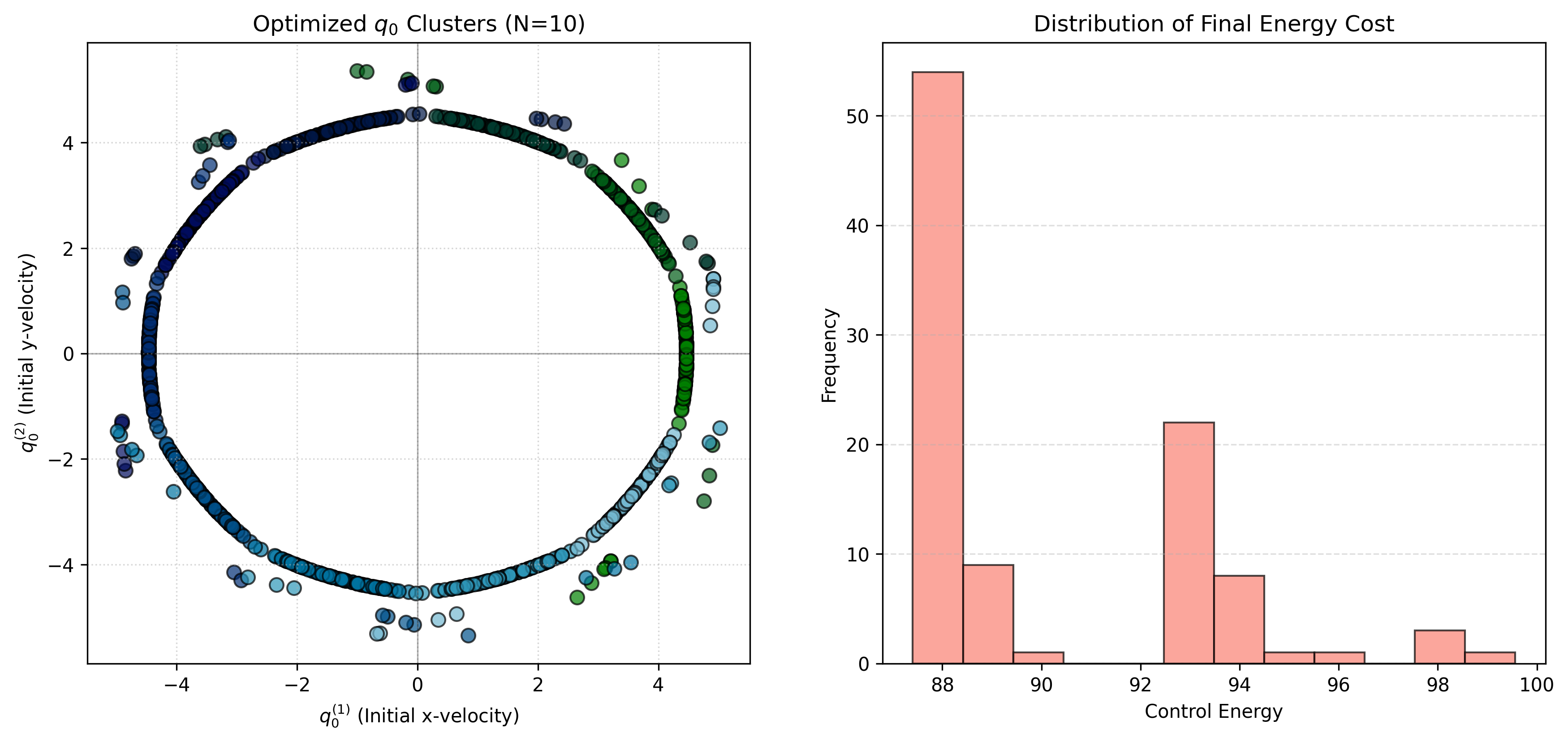}
        \caption{$N=10$, Ring Target}
        \label{fig:mc_n10_ring}
    \end{subfigure}
    \hfill
    \begin{subfigure}{0.48\textwidth}
        \centering
        \includegraphics[width=\linewidth]{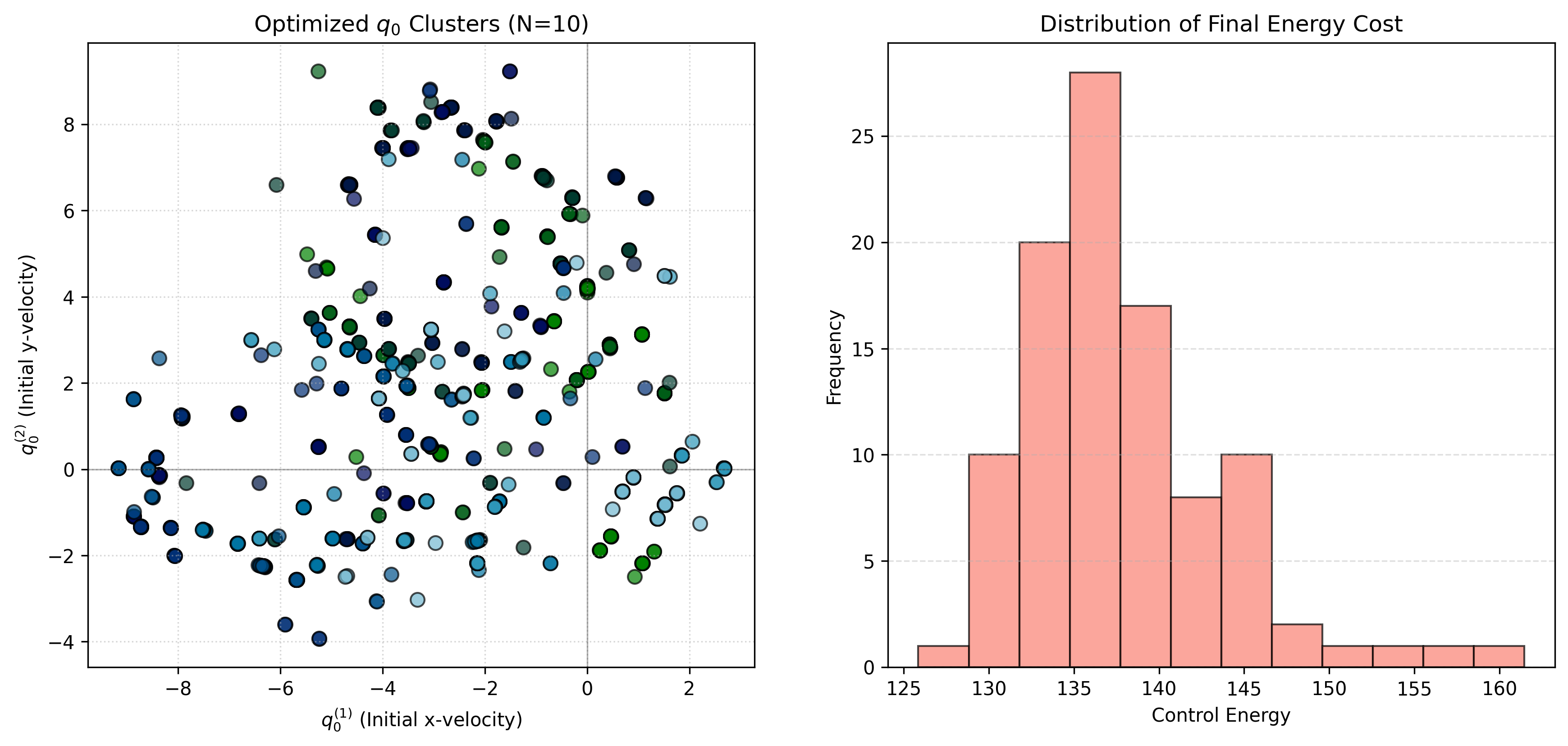}
        \caption{$N=10$, Heart Target}
        \label{fig:mc_n10_heart}
    \end{subfigure}

    \vspace{0.5cm} 

    \begin{subfigure}{0.48\textwidth}
        \centering
        \includegraphics[width=\linewidth]{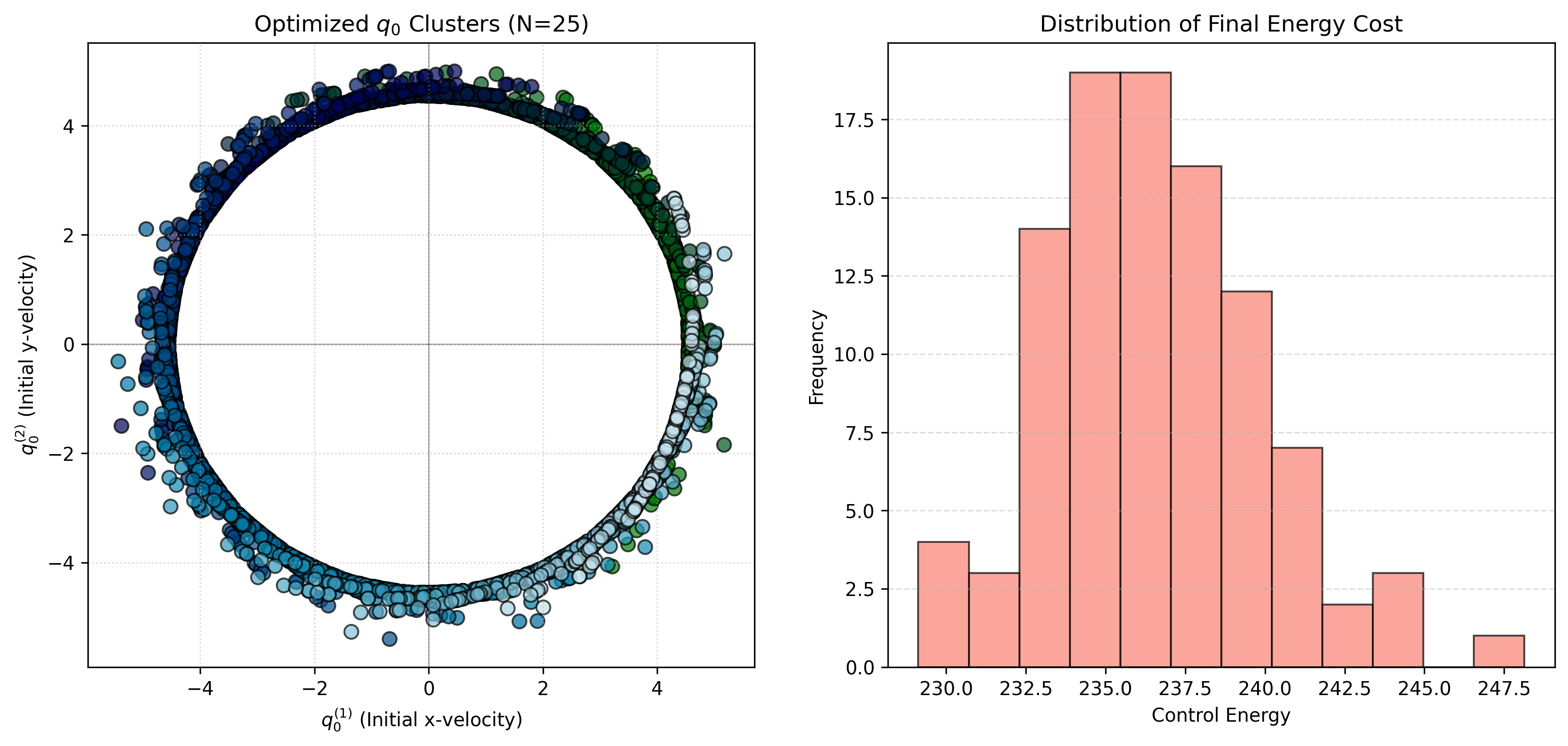}
        \caption{$N=25$, Ring Target}
        \label{fig:mc_n25_ring}
    \end{subfigure}
    \hfill
    \begin{subfigure}{0.48\textwidth}
        \centering
        \includegraphics[width=\linewidth]{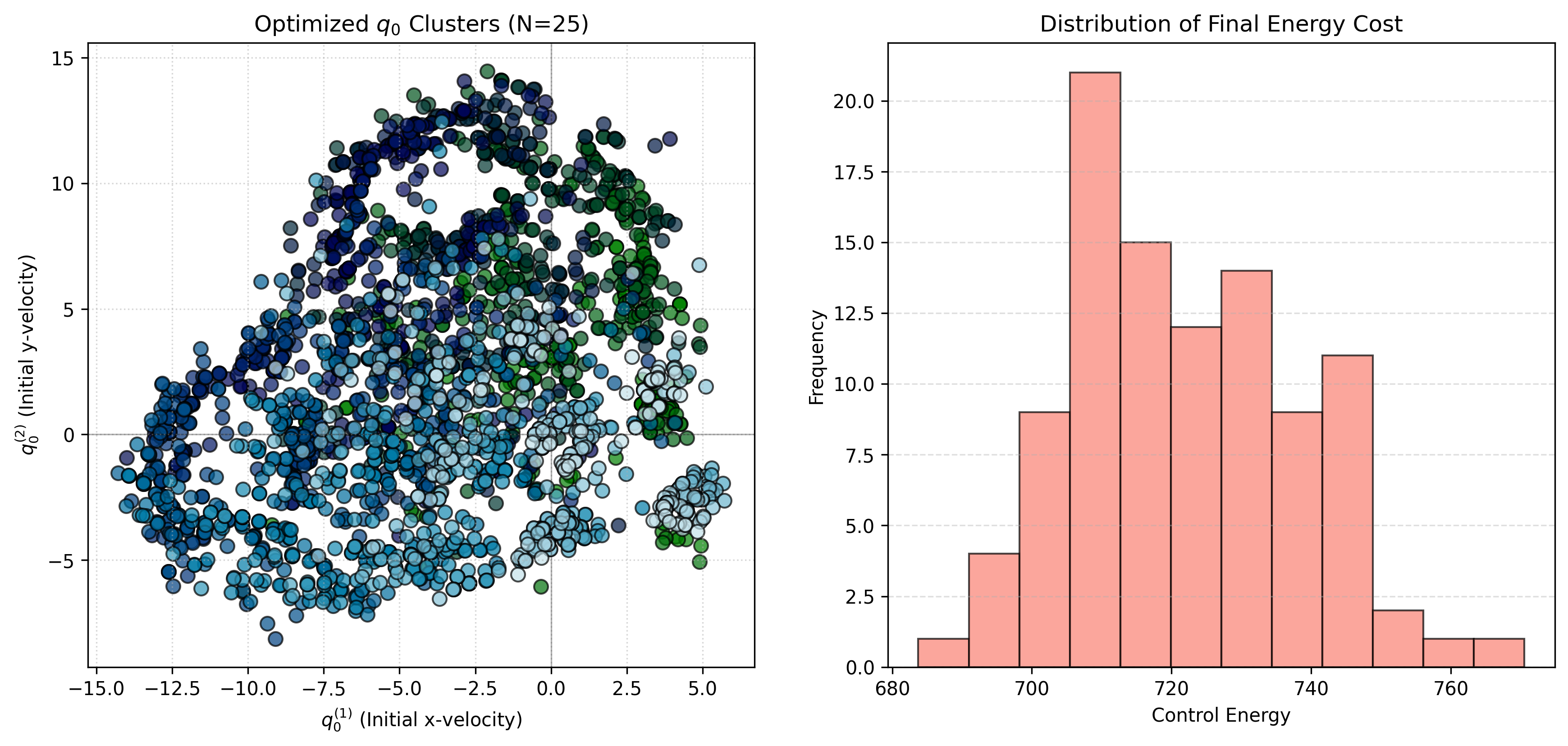}
        \caption{$N=25$, Heart Target}
        \label{fig:mc_n25_heart}
    \end{subfigure}

    \caption{Analysis of optimization landscape and control effort for $N \in \{10, 25\}$ agents within a point repeller background flow. 
    Each panel consists of: (Left) a scatter plot of optimized initial velocity $\mathbf{q}_{0,\text{opt}}$ across $M=100$ trials, where tight clustering per agent (color-coded) indicates convergence to consistent navigational strategies despite random initialization; 
    and (Right) the frequency distribution of final control energy $E_{\textrm{WHF}} $, indicating numerical stability across different target geometries. In plots (b) and (d), the rotated heart shape forms due to the background flow and the target distribution.}
    \label{fig:MC_analysis_full}
\end{figure}

  The left panels in Figure \ref{fig:MC_analysis_full} show the distribution of the optimized initial velocities $\mathbf{q}_{0,\text{opt}}$ across agents, visualized via scatter plots to identify emergent navigational strategies. The right panels show histograms of the energy expenditures, demonstrating the distribution of control energy across the random initializations.

\begin{figure}
    \centering
    \includegraphics[height=2.5in, trim={0 0 2.5cm 0}, clip]{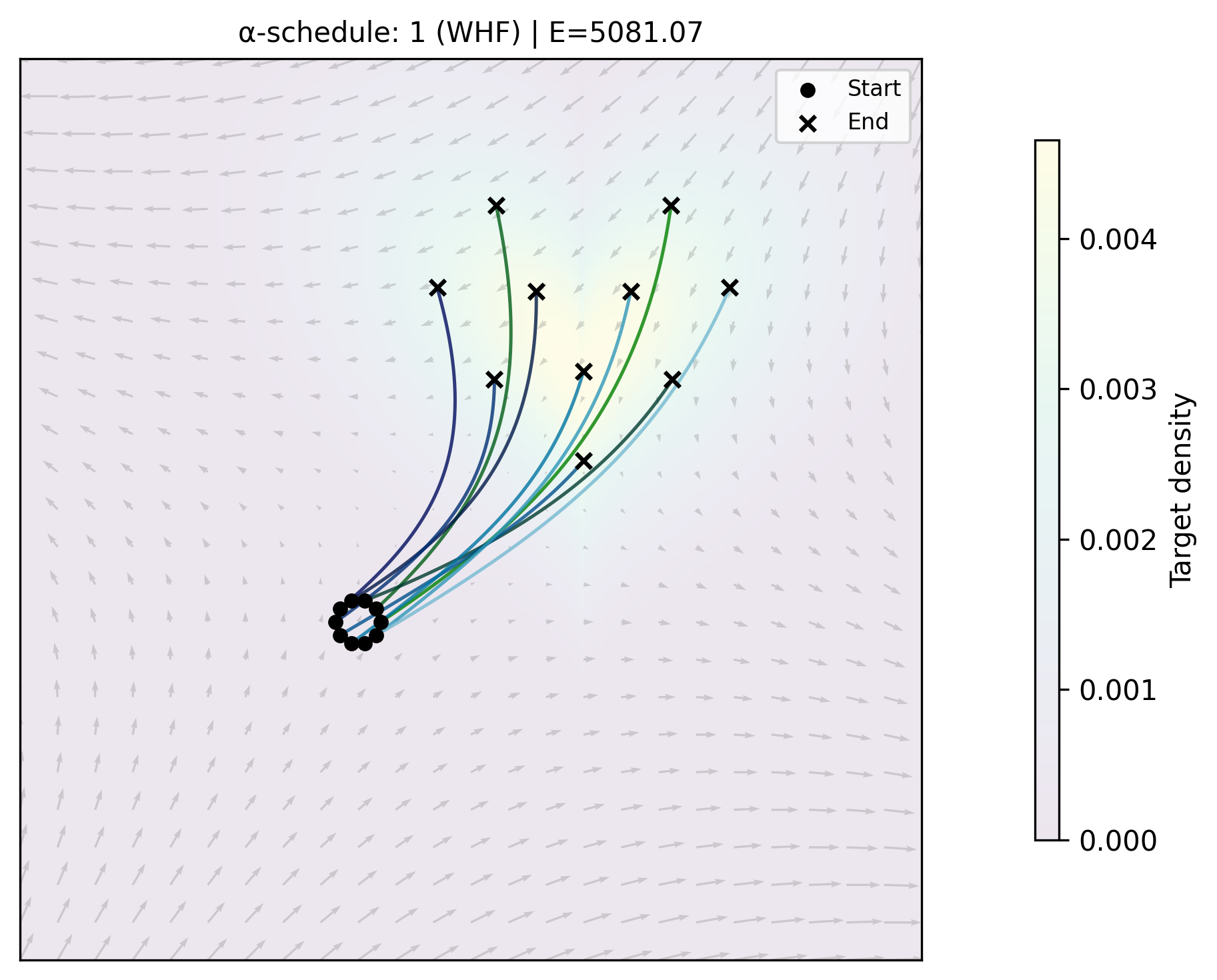}\includegraphics[height=2.5in]{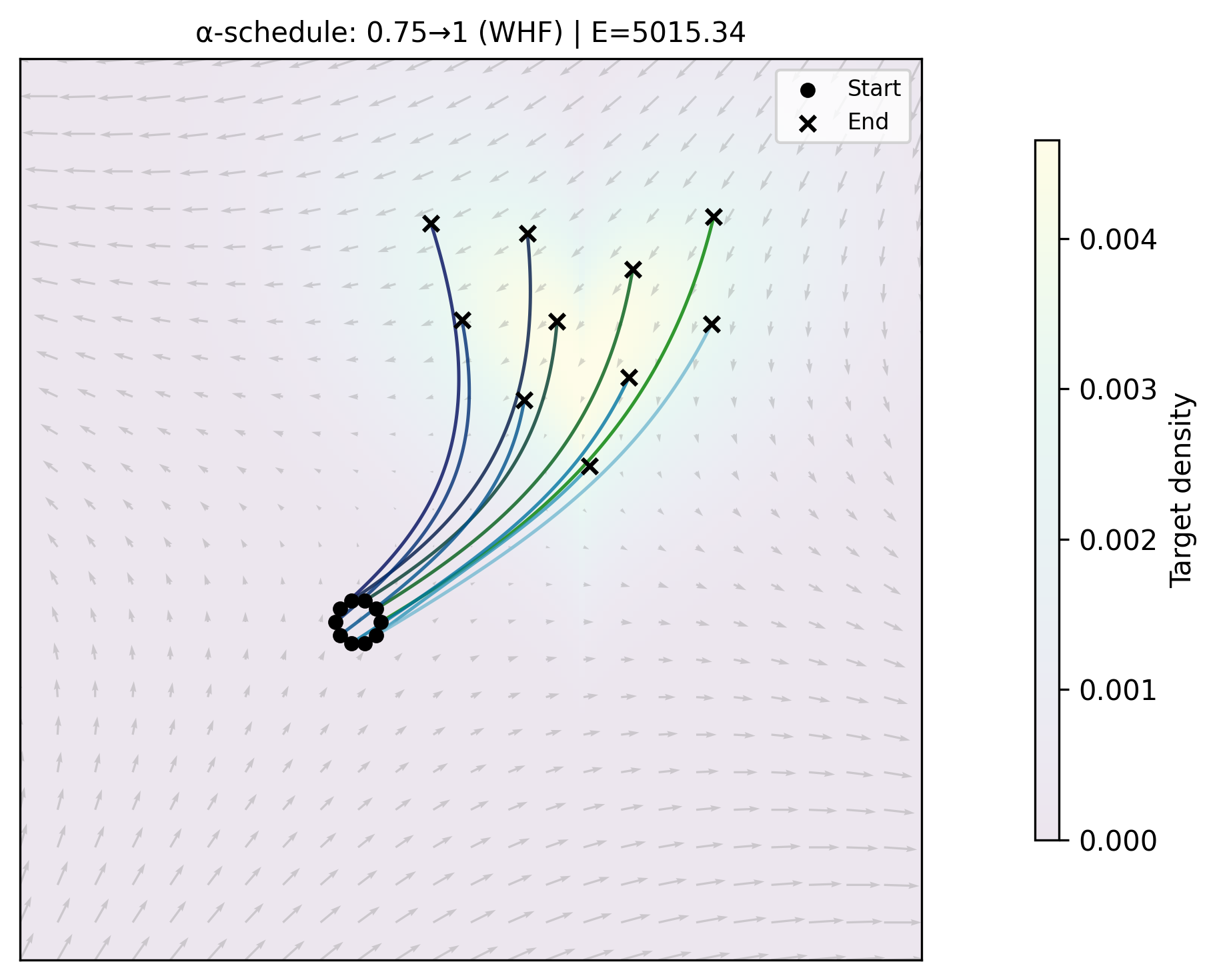}
    \caption{Applying a continuation (homotopy) method for multi-agent path-planning in a moving medium can result in further energy savings.}
    \label{fig:homotopy}
\end{figure}

One way to address local minima is to employ a homotopy-based continuation method \cite{Sun1995Homotopy-continuationEnergy, Richter1983ContinuationApplications}, which progressively deforms a simpler auxiliary problem into the target objective to ensure global convergence. Figure \ref{fig:homotopy} demonstrates how a continuation (homotopy) method for multi-agent path-planning in a moving medium can result in further energy savings. 
The objective functional $\mathcal{J}$ is generally non-convex due to the nonlinear term $w$. To avoid local minima, we employed a parameter continuation strategy (homotopy) using a family of flow fields $w_\alpha(t, x) = \alpha w(t, x)$ for $\alpha \in [0, 1]$. Further investigation of the adaptive parameters in homotopy methods is needed to guarantee improved performance. The left plot in Figure \ref{fig:homotopy} demonstrates the outcome using a singleton family $\alpha\in \{1\}$, and the right plot demonstrates the outcome using the family $\alpha\in \{.75,1\}$, resulting in energy savings.

\section{Future work}\label{sec:future}

Based on a modified optimal transport formulation and its WHF, we designed an optimization strategy with ODE constraints for multi-agent path-planning in a moving medium. Compared with the existing trajectory-based methods, this modified model features a finite dimensional search space, and the problem is solved by a combination of existing optimization methods and ODE schemes. Numerical experiments demonstrated its performance in various scenarios. Meanwhile, the model invites several questions that are of interest for further investigation. 

From a theoretical perspective, the feasibility, namely the existence of solutions for the WHF to achieve the target distribution, is a critical issue for the model. Our numerical experience shows that the model is robust in many situations. However, it is possible to set up pathological flow fields or desirable distributions so that the model becomes infeasible. Identifying the conditions that describe the feasibility is important and difficult. Furthermore, understanding the theoretical connections between this model and existing ones is also interesting. 

From a practical point of view, there are a number of considerations. For example, collision avoidance and decentralization in the planning. Collision avoidance includes avoiding obstacles or other agents. Obstacles can be incorporated through a potential functional in the formulation as being done in many optimal transport and mean field game problems.  Collision avoidance with other agents is theoretically guaranteed by the uniqueness of ODEs, which ensures that the trajectories do not intersect with each other if their initial locations are different and the right hand side of the WHF is Lipschitz continuous. However, this is only the case if the agents are abstracted as points. In reality, agents have given sizes. This requires building buffer zones around them, which may be achieved using interactive potentials to repell agents when they get close to each other in the formulation. The current model and methods provide the path-planning in advance, and they do not support the computation on the fly. For this reason, decentralized strategies are more helpful. 

In computation, the model is finite dimensional, which has the potential to lead to more efficient computation. Combinations of the of-the-shelf optimization methods and ODE solvers, as we present in this paper, are not necessarily efficient. Specialized algorithms that take advantage of the properties of the model are desirable to further reduce the computational cost, scale up in the number of agents, and improve stability. Our experiments indicate that the global optimizers are achieved in some examples, such as the ``ring'' target distribution, but it is more likely that the algorithm only obtains local optimizers in other examples. Hence, developing algorithms and theory to obtain the global solutions remains as a difficult task.

\bibliographystyle{abbrv}
\bibliography{references}

\end{document}